\documentclass[mathpazo]{1586312043190378}

\usepackage{amsmath}
\usepackage{enumerate}
\usepackage{mathrsfs}
\usepackage{color}
\usepackage{tikz}
\usepackage{mathtools}
\usepackage{breqn}
\usepackage{colortbl}
\usepackage{subfigure}
\usepackage{xr}

\newcommand{\defword}[1]{\textcolor{red}{\em #1}}

\begin{document}
	\title{Multistability of Bi-Reaction Networks}
	
	\author[O.~Author]{Yixuan Liang, Xiaoxian Tang and Qian Zhang}
	\address{School of Mathematical Sciences, Beihang University,
		Beijing 100191, P.R. China.\\
		Key Laboratory of Mathematics, Informatics and Behavioral Semantics, Ministry of Education,  Beijing 100191, P. R. China.\\
		State Key Laboratory of Software Development Environment, Beihang University, Beijing 100083, P.R. China.}



	\begin{abstract}	
		We provide a sufficient and necessary condition in terms of the stoichiometric coefficients for a bi-reaction network to admit multistability. Also, this result 
  completely characterizes  the bi-reaction networks according to if 
  they admit multistability. 
	\end{abstract}
	
	
	\maketitle

	\section{Introduction}\label{intro}
	This work addresses the multistability problem for the dynamical systems arising from the bio-chemical reaction networks (under mass-action kinetics). The problem is 
 how to efficiently determine if a reaction  network admits at least two stable positive steady states in the same stoichiometric compatibility class.  Multistability  is important in mathematical biology since it widely exists in
the decision-making process and switch-like behavior in cellular signaling (e.g., \cite{BF2001, XF2003, CTF2006, SN2020,SH2021}). 
In practice, one way to detect multistability is to first find nondegenerate multistationarity (i.e., to check if the network admits more than one positive nondegenerate steady state). 
Usually, one can  obtain two stable steady states if the number of positive nondegenerate steady states  is at least three  (e.g., \cite{OSTT2019,TF2020,CH2021}). Generally,  
deciding multistationarity/multistability or computing the witnesses (i.e., a choice of parameters for which the network exhibits multistationarity/multistability)  is challenging because the problem is known to be a special real quantifier elimination problem (that means we want to efficiently obtain the information of real solutions of a semi-algebraic system, e.g., \cite{HTX2015, B2017}). However, 
 there indeed exists a collection of efficient/practical methods for  detecting multistationarity  (e.g., \cite{ShinarFeinberg2012,CF2012,  signs, joshi2013complete}).  Most of these approaches are  to check if the determinant of a certain Jacobian matrix changes sign \cite{CF2005, BP,WiufFeliu_powerlaw,CFMW, DMST2019,SadeghimaneshFeliu2019}. 

 One big goal in the area of reaction network  is to look for the ``explicit" criteria. That means we hope to tell the dynamical behaviors of a network by reading the network itself without doing any expensive  computations.  One typical result, which makes the big goal realistic, is the well-known deficiency zero theorem and the deficiency one theorem \cite{1.24 add_1}. 
 So far, such explicit criteria for detecting   multistationarity/multistability are only known for small networks with one species or up to two reactions (possibly reversible) \cite{Joshi:Shiu:Multistationary, shiu-dewolff}.
For instance, in \cite{Joshi:Shiu:Multistationary},  the authors completely characterized one-species networks by ``arrow diagrams", and the number of (stable) steady states can be read off by looking at the existence of $T$-alternating subnetworks with certain type of arrow diagrams. Later, in \cite{the first}, the criterion for multistationarity described by  arrow diagrams is extended to more general networks with one-dimensional stoichiometric subspaces.    
Since for the one-dimensional networks,  admitting at least three positive steady states is a necessary condition for admitting multistability (e.g., \cite[Theorem 3.4]{txzs2020}),  the explicit criterion for admitting three positive steady states (described by  ``bi-arrow diagrams")  is studied in  \cite{linkexin}. Also, the authors of 
\cite{linkexin} has completely characterized the stoichiometric coefficients of the bi-reaction networks that admit at least three positive steady states. We remark that in the point of view of real algebraic geometry, an explicit criterion for multistationarity/multistability  is essentially an explicit criterion for deciding number of real solutions of a special class of semi-algebraic systems. Some related recent work is the extension of the Descartes' rule of signs for the high dimensional algebraic systems (e.g., \cite{1.24 add_2, 1.24 add_3}), which can also be applied to the steady-state systems arising from bio-chemical reaction networks. 

In this paper, we focus on the bi-reaction networks that admit finitely many positive steady states. The main result is an explicit criterion for deciding multistability of bi-reaction networks (see Theorem \ref{thm:main}). By this result, we completely classify all non-trivial bi-reaction networks according to if they admit multistability or not (here, a ``non-trivial" network means this network admits at least one positive steady state). 
This work can be viewed as an extension of 
\cite{linkexin}, since in \cite{linkexin}, all bi-reaction networks are classified according to if they admit at least three positive steady states (recall that admitting three positive steady states is a necessary condition for multistability).  

The rest of this paper is organized as follows.
In Section~\ref{sec:back}, we recall the basic definitions and notions for the
reaction networks and the multistationarity/multistability. 
In Section \ref{sec:main}, we present the main theorem (a sufficient and necessary condition in terms of the stoichiometric coefficients for a bi-reaction network to admit multistability), and we illustrate how to use the theorem for deciding   multistability by several examples.  In Section \ref{sec:proof}, we present the proof of the main theorem by discussing several cases.  In the supplementary materials \footnote{https://github.com/65536-1024/one-dim}, we present a list of useful lemmas and their proofs.

	\section{Background}\label{sec:back}
	
	\subsection{Chemical reaction networks}\label{sec:pre}
	
	
	
	In this paper, we follow the standard notions on reaction networks used in \cite{linkexin, txzs2020}.
	A \defword{reaction network}  $G$  (or \defword{network} for short) consists of finitely many reactions:
	\begin{align}\label{eq:network}
		\alpha_{1j}X_1 +
		\dots +
		\alpha_{sj}X_s
		~ \xrightarrow{\kappa_j} ~
		\beta_{1j}X_1 +
		\dots +
		\beta_{sj}X_s,
		\;
		\;\;\;\;~
		j=1,2, \ldots, m,
	\end{align}
	where $X_1, \ldots, X_s$ denote $s$ \defword{species},  the \defword{stoichiometric coefficients} $\alpha_{ij}$ and $\beta_{ij}$ are non-negative integers, each $\kappa_j \in \mathbb R_{>0}$ is a \defword{rate constant} corresponding to the
	$j$-th reaction,
	and
	we assume that
	\begin{align}\label{eq:netcon}
		\text{for every}\; j\in \{1, \cdots, m\},\; (\alpha_{1j},\cdots,\alpha_{sj})\neq (\beta_{1j},\cdots,\beta_{sj}).
	\end{align}
	The
	\defword{stoichiometric matrix} of
	$G$,
	denoted by ${\mathcal N}$, is the
	$s\times m$ matrix with
	$(i, j)$-entry equal to $\beta_{ij}-\alpha_{ij}$.
	The \defword{stoichiometric subspace}, denoted by $S$, is the real vector space spanned by the column vectors of ${\mathcal N}$.

	The concentrations of the species $X_1,X_2, \ldots, X_s$ are denoted by $x_1, x_2, \ldots, x_s$, respectively. Note that
	$x_i$ can be considered as a function in the time variable $t$.
	Under the assumption of mass-action kinetics, we describe how these concentrations change  in $t$ by the following system of ordinary differential equations
	(ODEs):
	\begin{align}\label{eq:sys}
		\dot{x}~=~(f_1(\kappa; x), \cdots, f_s(\kappa; x))^{\top}~:=~{\mathcal N}\cdot \begin{pmatrix}
			\kappa_1 \, \prod \limits_{i=1}^s x_i^{\alpha_{i1}} \\
			\kappa_2 \, \prod \limits_{i=1}^s x_i^{\alpha_{i2}} \\
			\vdots \\
			\kappa_m \, \prod \limits_{i=1}^s x_i^{\alpha_{im}} \\
		\end{pmatrix}~,
	\end{align}
	where $x$ denotes the vector $(x_1, x_2, \ldots, x_s)$,
	and $\kappa$ denotes the vector $(\kappa_1, \ldots, \kappa_m)$.
	Note that for every $i\in \{1, \ldots, s\}$, $f_{i}(\kappa;x)$ is a polynomial  in $\mathbb Q[\kappa, x]$.
	
	A \defword{conservation-law matrix} of $G$, denoted by $W$, is any row-reduced $d\times s$ matrix (here, $d:=s-{\rm rank}({\mathcal N})$), whose rows form a basis of $S^{\perp}$. Note that ${\rm rank}(W)=d$. Especially, if the stoichiometric subspace of $G$ is one-dimensional, then
	${\rm rank}({\mathcal N})=1$ and ${\rm rank}(W)=s-1$.
	Note that the system~\eqref{eq:sys} satisfies $W \dot x =0$,  and 
	any trajectory $x(t)$ beginning at a nonnegative vector $x(0)=x^0 \in
	\mathbb{R}^s_{> 0}$ remains, for all positive time,
	in the following \defword{stoichiometric compatibility class} with respect to the  \defword{total-constant vector} $c:= W x^0 \in {\mathbb R}^d$:
	\begin{align}\label{eq:pc}
		{\mathcal P}_c~:=~ \{x\in {\mathbb R}_{\geq 0}^s \; :\; Wx=c\}.~
	\end{align}
	
	\subsection{Multistationarity and multistability}\label{sec:mm}
	For a given rate-constant vector $\kappa\in \mathbb{R}_{>0}^m$, a \defword{steady state} 
	of~\eqref{eq:sys} is a concentration vector
	$x^* \in \mathbb{R}_{\geq 0}^s$ such that $f_1(\kappa, x^*)=\cdots=f_s(\kappa, x^*)=0$, where  $f_1, \ldots, f_s$ are on the
	right-hand side of the
	ODEs~\eqref{eq:sys}.
	If all coordinates of a steady state $x^*$ are strictly positive (i.e., $x^*\in \mathbb{R}_{> 0}^s$), then we call $x^*$ a \defword{positive steady state}. 
	We say a steady state $x^*$ is \defword{nondegenerate} if
	${\rm im}\left({\rm Jac}_f (x^*)|_{S}\right)=S$,
	where ${\rm Jac}_f(x^*)$ denotes the Jacobian matrix of $f$, with respect to $x$, at $x^*$.
	A steady state $x^*$ is \defword{exponentially stable} (or simply \defword{stable})
	if it is nondegenerate, and  all non-zero eigenvalues of ${\rm Jac}_f(x^*)$ have negative real parts.
	Note that if a steady state is exponentially stable, then it is locally asymptotically stable \cite{P2001}.
	
	
	
	
	Suppose $N\in {\mathbb Z}_{\geq 0}$. We say a  network  \defword{admits $N$ (nondegenerate) positive steady states}  if there exist a rate-constant
	vector $\kappa\in  \mathbb{R}_{>0}^m$ and a total-constant vector $c\in  \mathbb{R}^d$ such that it has $N$ (nondegenerate) positive steady states  in the stoichiometric compatibility class ${\mathcal P}_c$.
	Similarly, we say a  network  \defword{admits $N$ stable positive steady states}  if there exist a rate-constant
	vector $\kappa\in  \mathbb{R}_{>0}^m$ and a total-constant vector $c\in  \mathbb{R}^d$ such that  it has $N$ stable positive steady states  in ${\mathcal P}_c$.
	
	The \defword{maximum number of positive steady states} of a network $G$ is
	{\footnotesize
		\[cap_{pos}(G)\;:=\;\max\{N\in {\mathbb Z}_{\geq 0}\cup \{+\infty\}\; :\; G \;\text{admits}\; N\; \text{positive steady states}\}.\]
	}Similarly, we define
	{\footnotesize
		\[cap_{nondeg}(G)\;:=\;\max\{N\in {\mathbb Z}_{\geq 0}\cup \{+\infty\}\; :\; G \;\text{admits}\; N\; \text{nondegenerate positive steady states}\}\]
	}and
	{\footnotesize
		\[cap_{stab}(G)\;:=\;\max\{N\in {\mathbb Z}_{\geq 0}\cup \{+\infty\}\; :\; G \;\text{admits}\; N\; \text{stable positive steady states}\}.\]
	}We say a network admits \defword{multistationarity} if $cap_{pos}(G)\geq 2$.
	We say a network admits \defword{nondegenerate multistationarity} if $cap_{nondeg}(G)\geq 2$.
	We say a network admits \defword{multistability} if $cap_{stab}(G)\geq 2$.
	\begin{theorem}\label{thm:nc}\cite[Theorem 6.1]{txzs2020}
		Given a network $G$ with a one-dimensional stoichiometric subspace,
		if $cap_{pos}(G)<+\infty$, then $cap_{nondeg}(G)=cap_{pos}(G)$.
	\end{theorem}
	
	

	\section{Main result}\label{sec:main}
	In this section, we focus on the bi-reaction network $G$: 
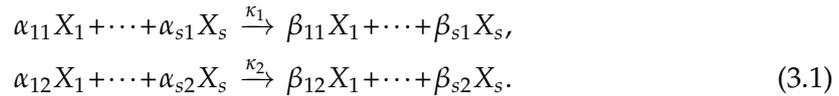
	\begin{align}\label{eq:2network}
		\alpha_{11}X_1 +
		\dots +
		\alpha_{s1}X_s
		&~ \xrightarrow{\kappa_1} ~
		\beta_{11}X_1 +
		\dots +
		\beta_{s1}X_s,
		\notag \\
		\alpha_{12}X_1 +
		\dots +
		\alpha_{s2}X_s
		&~ \xrightarrow{\kappa_2} ~
		\beta_{12}X_1 +
		\dots +
		\beta_{s2}X_s.
	\end{align}
	First, we define the following sets of indices, which give a partition of the set $\{1, \ldots, s\}$:
	\begin{align}\label{eq:s}
		S_1&:=\{i\;:\;\alpha_{i1}>\alpha_{i2},\ \beta_{i1}>\alpha_{i1},\ 1\leqslant i\leqslant s \},\notag \\
		S_2&:=\{i\;:\;\alpha_{i1}<\alpha_{i2},\ \beta_{i1}<\alpha_{i1},\ 1\leqslant i\leqslant s \},\notag \\
		S_3&:=\{i\;:\;\alpha_{i1}>\alpha_{i2},\ \beta_{i1}<\alpha_{i1},\ 1\leqslant i\leqslant s \},\notag \\
		S_4&:=\{i\;:\;\alpha_{i1}<\alpha_{i2},\ \beta_{i1}>\alpha_{i1},\ 1\leqslant i\leqslant s \},\;\;\;\;\text{and}\notag \\
		S_5&:=\{i\;:\;\alpha_{i1}=\alpha_{i2}\ \text{or}\ \beta_{i1}=\alpha_{i1},\ 1\leqslant i\leqslant s \}.
	\end{align}
	For each $i\in \{1, \ldots, s\}$, we define the following notions:
	\begin{align}\label{eq:ar}
		a_i:=|\alpha_{i1}-\alpha_{i2}|,\; \;\;
		\gamma_i:=|\beta_{i1}-\alpha_{i1}|.      
	\end{align}

	\begin{theorem}\label{thm:main}
		Given a bi-reaction network $G$ \eqref{eq:2network} with a one-dimensional stoichiometric subspace,  suppose $0<cap_{pos}(G)<+\infty$.
		\begin{itemize}
			\item[(a)] If all the four sets $S_1,\; S_2,\; S_3,\; S_4$ are non-empty, then $G$ admits  multistability if and only if 
			\begin{align}\label{end4}
				\sum\limits_{i\in S_1}a_i > \min\limits_{i\in S_4} \{a_i \},\; or \; \sum\limits_{i\in S_2}a_i > \min\limits_{i\in S_3} \{a_i \}.
			\end{align}
			
			\item[(b)] If there are exactly three of the four sets $S_1,\; S_2,\; S_3,\; S_4$ are non-empty, then $G$ admits  multistability if and only if one of the following four statements (1)--(4) holds.
			
			(1) $S_1$, $S_3$ and $S_4$ are non-empty, and 
			\begin{align}
				\sum\limits_{i\in S_1}a_i>\min\limits_{i\in S_{4}}\{a_i\}.
			\end{align}
			
			(2) $S_2$, $S_3$ and $S_4$ are non-empty, and 
			\begin{align}
				\sum\limits_{i\in S_2}a_i>\min\limits_{i\in S_{3}}\{a_i\}.
			\end{align}
			
			(3) $S_1$, $S_2$ and $S_3$ are non-empty, and there exists a subset $S_2^{*}$ of $S_2$ such that
			\begin{align}
				\sum\limits_{i\in S_{3}}a_i>\sum\limits_{i\in S_2^{*}}a_i>\min\limits_{i\in S_{3}}\{a_i\}.
			\end{align}
			
			(4) 
			$S_1$, $S_2$ and $S_4$ are non-empty, and there exists a subset $S_1^{*}$ of $S_1$ such that
			\begin{align}
				\sum\limits_{i\in S_{4}}a_i>\sum\limits_{i\in S_1^{*}}a_i>\min\limits_{i\in S_{4}}\{a_i\}.
			\end{align}
			
			\item[(c)] If there are exactly two of the four sets $S_1,\; S_2,\; S_3,\; S_4$ are non-empty, then $G$ admits  multistability if and only if one of the following two statements (1)--(2) holds.
			(1) $S_2$ and $S_3$ are non-empty and there exists a subset $S_2^{*}$ of $S_2$ such that
			\begin{align}
				\sum\limits_{i\in S_{3}}a_i>\sum\limits_{i\in S_2^{*}}a_i>\min\limits_{i\in S_{3}}\{a_i\}.
			\end{align}
			
			(2) $S_1$ and $S_4$ are non-empty, and there exists a subset $S_1^{*}$ of $S_1$ such that
			\begin{align}
				\sum\limits_{i\in S_{4}}a_i>\sum\limits_{i\in S_1^{*}}a_i>\min\limits_{i\in S_{4}}\{a_i\}.
			\end{align}

			\item[(d)] If only one of the four sets $S_1,\; S_2,\; S_3,\; S_4$ is non-empty, then $G$ admits no  multistability.
		\end{itemize}
	\end{theorem}

\begin{example}The following examples illustrate how Theorem \ref{thm:main} works.
  \begin{itemize}
    \item [(a)] Consider the following network:$$G:\begin{aligned}
      4X_1+X_2+X_3&\rightarrow 5X_1+X_4\\
      X_1+2X_2+X_4&\rightarrow 3X_2+X_3
    \end{aligned}$$
It is straightforward to check that:
    $$\begin{aligned}
      &a_1=3,~ a_2=1, ~a_3=1,~a_4=1.\\
      &S_1=\{1\},~ S_2=\{2\},~ S_3=\{3\}, ~ S_4=\{4\},~ S_5=\emptyset.\\
      &\sum_{i\in S_1} a_i=a_1=3>1=a_4=\min_{i\in S_4}\{a_i\}.
    \end{aligned}$$ So by Theorem \ref{thm:main} (a), we have $cap_{stab}(G)\geq 2$. 
     The steady-state system augmented with the conservation laws is
     \begin{align*}
\kappa_1 x_1^4 x_2 x_3-\kappa_2 x_1 x_2^2 x_4=0,\\
-x_1-x_2-c_1=0,\\
-x_1-x_3-c_2=0,\\
x_1-x_4-c_3=0.
     \end{align*}
    One can check that for $c_1=-2,~c_2=-\dfrac{17}{10},~c_3=\dfrac{3}{10},~\kappa_1=1,$ and  $\kappa_2=1$,
    the network has $3$ positive steady states: $$\begin{aligned}
      &x^{(1)}=(0.3293,~1.671,~1.371,~0.02930),\\
      &x^{(2)}=(1.000,~1.000,~0.7000,~0.7000),\\
      &x^{(3)}=(1.548,~0.4521,~0.1521,~1.248),
    \end{aligned}$$
    where $x^{(1)}$ and $x^{(3)}$ are stable.
		\item [(b) (i)] Consider the following network:$$G:\begin{aligned}
			3X_1+3X_2+2X_3+X_4&\rightarrow 4X_1+4X_2+X_3+2X_4\\
			X_1+X_2+X_3+4X_4&\rightarrow 2X_3+3X_4
		\end{aligned}$$
		It is straightforward to check that:
		$$\begin{aligned}
			&a_1=2,~ a_2=2, ~a_3=1, ~a_4=3.\\
			&S_1=\{1,\;2\},~ S_2=\emptyset,~ S_3=\{3\}, ~ S_4=\{4\},~ S_5=\emptyset.\\
			&\sum_{i\in S_1} a_i=a_1+a_2=4>3=a_4=\min_{i\in S_4}\{a_i\}.
		\end{aligned}$$ 
  So by Theorem \ref{thm:main} (b) (1), we have $cap_{stab}(G)\geq 2$. 
     The steady-state system augmented with the conservation laws is
     \begin{align*}
\kappa_1 x_1^3 x_2^3 x_3^2 x_4-\kappa_2 x_1 x_2 x_3 x_4^4=0,\\
x_1-x_2-c_1=0,\\
x_1+x_3-c_2=0,\\
x_1-x_4-c_3=0.
\end{align*}
		One can check that for                           $c_1=\dfrac{9}{100},~c_2=3,~c_3=\dfrac{1}{10},~\kappa_1=1,$ and  $\kappa_2=2$,
		the network has $3$ positive steady states:  $$\begin{aligned}
			&x^{(1)}=(0.1448,~0.05478,~2.855,~0.04478),\\
   &x^{(2)}=(0.7442,~0.6542,~2.256,~0.6442),\\
			&x^{(3)}=(2.103,~2.013,~0.8967,~2.003),
		\end{aligned}$$where $x^{(1)}$ and $x^{(3)}$ are stable.
	\item [(b) (ii)] Consider the following network:$$G:\begin{aligned}
		2X_1+2X_2+X_3+X_4+X_5+3X_6&\rightarrow 3X_1+X_2\\
		X_1+3X_2+3X_3+2X_4&\rightarrow 4X_2+4X_3+3X_4+X_5+3X_6
	\end{aligned}$$
	It is straightforward to check that:
	$$\begin{aligned}
		&a_1=1,~ a_2=1, ~a_3=2, ~a_4=1,~a_5=1,~a_6=3.\\
		&S_1=\{1\},~ S_2=\{2,\;3,\;4\},~ S_3=\{5,\;6\}, ~ S_4=\emptyset,~ S_5=\emptyset.
	\end{aligned}$$
Choose $S_2^*=\{2,\;3\}$. Note that
$$\begin{aligned}
	\sum\limits_{i\in S_{3}}a_i=4>\sum\limits_{i\in S_2^{*}}a_i=3>\min\limits_{i\in S_{3}}\{a_i\}=1.
\end{aligned}$$
So by Theorem \ref{thm:main} (b) (3), we have $cap_{stab}(G)\geq 2$. 
     The steady-state system augmented with the conservation laws is
     \begin{align*}
\kappa_1 x_1^2 x_2^2 x_3 x_4 x_5 x_6^3-\kappa_2 x_1 x_2^3 x_3^3 x_4^2=0,\\
x_1+x_2-c_1=0,\\
x_1+x_3-c_2=0,\\
x_1+x_4-c_3=0,\\
x_1+x_5-c_4=0,\\
3x_1+x_6-c_5=0.
\end{align*}
		One can check that for                           $c_1=101,~c_2=101,~c_3=1000,~c_4=100,~c_5=315,~\kappa_1=1,$ and  $\kappa_2=72$,
		the network has $4$ positive steady states:  $$\begin{aligned}
			&x^{(1)}=(32.09,~68.91,~68.91,~967.9,~67.91,~218.7),\\
   &x^{(2)}=(86.24,~14.76,~14.76,~913.8,~13.76,~56.29),\\
			&x^{(3)}=(97.55,~3.450,~3.450,~902.5,~2.450,~22.35),\\
&x^{(4)}=(99.54,~1.464,~1.464,~900.5,~0.4641,~16.39),
		\end{aligned}$$where $x^{(2)}$ and $x^{(4)}$ are stable.
\item [(c)] Consider the following network:$$G:\begin{aligned}
	3X_1+4X_2+5X_3+2X_4&\rightarrow 4X_1+5X_2+7X_3+3X_4+X_5\\
	2X_1+2X_2+4X_3+3X_4+2X_5&\rightarrow X_4
\end{aligned}$$
It is straightforward to check that:
$$\begin{aligned}
	&a_1=1,~ a_2=2, ~a_3=1, ~a_4=1,~a_5=2.\\
	&S_1=\{1,\;2,\;3\},~ S_2=\emptyset,~ S_3=\emptyset, ~ S_4=\{4,\;5\},~ S_5=\emptyset.
\end{aligned}$$
Choose $S_1^*=\{2\}$. Note that
$$\begin{aligned}
	\sum\limits_{i\in S_{4}}a_i=3>\sum\limits_{i\in S_1^{*}}a_i=2>\min\limits_{i\in S_{4}}\{a_i\}=1.
\end{aligned}$$
So by Theorem \ref{thm:main} (c) (1), we have $cap_{stab}(G)\geq 2$. 
     The steady-state system augmented with the conservation laws is
     \begin{align*}
\kappa_1 x_1^3 x_2^4 x_3^5 x_4^2-2\kappa_2 x_1^2 x_2^2 x_3^4 x_4^3 x_5^2=0,\\
x_1-x_2-c_1=0,\\
2x_1-x_3-c_2=0,\\
x_1-x_4-c_3=0,\\
x_1-x_5-c_4=0.\\
\end{align*}
		One can check that for                           $c_1=100,~c_2=1,~c_3=101,~c_4=90,~\kappa_1=1,$ and  $\kappa_2=328$,
		the network has $4$ positive steady states:  $$\begin{aligned}
			&x^{(1)}=(101.6,~1.588,~202.2,~0.5879,~11.59),\\
   &x^{(2)}=(108.1,~8.081,~215.2,~7.081,~18.08),\\
			&x^{(3)}=(128.2,~28.21,~255.4,~27.21,~38.21),\\
&x^{(4)}=(190.6,~90.62,~380.2,~89.62,~100.6),
		\end{aligned}$$where $x^{(1)}$ and $x^{(3)}$ are stable.

  \end{itemize}
\end{example}

	\section{Proofs}\label{sec:proof}
	\newtheorem{assumption}{Assumption}[section]
	\begin{assumption}\label{assumption1}
		Without loss of generality, for any network $G$ defined in  \eqref{eq:2network},
		we  assume that the set  $S_5$ of indices  defined in \eqref{eq:s} is empty throughout the rest of the paper. In fact, 
  if $0<cap_{pos}(G)<+\infty$, one can always construct a new network such that the original network $G$ is dynamically equivalent to the new one and the set $S_5$ for the new network is empty (e.g., \cite[Lemma 5.1, Lemma 5.2]{txzs2020}). 
	\end{assumption}

	We assume that any bi-reaction network $G$ mentioned has the form \eqref{eq:2network}. 
 Notice that by Assumption \ref{assumption1}, we have \begin{align}\label{empty_1}
		\beta_{i1}-\alpha_{i1}\neq 0,\; \text{for}\;\; i=1,\cdots,s.
	\end{align}
If $G$ has a one-dimensional stoichiometric subspace, then
  there exists $\lambda\in \mathbb{R}\ (\lambda \neq 0)$ such that
	\begin{align}\label{reldef 1}
		\left(
		\begin{array}{c}
			\beta_{12}-\alpha_{12}\\
			\vdots\\
			\beta_{s2}-\alpha_{s2}
		\end{array}
		\right)\;=\;
		\lambda \left(
		\begin{array}{c}
			\beta_{11}-\alpha_{11}\\
			\vdots\\
			\beta_{s1}-\alpha_{s1}
		\end{array}
		\right).
	\end{align}
	By \cite[Lemma 4.1] {Joshi:Shiu:Multistationary} (also, see \cite[Lemma 4.2]{txzs2020}), we assume that $\lambda<0$ in \eqref{reldef 1} (otherwise, the network admits no positive steady state).
	By substituting
	\eqref{reldef 1} into $f_1, \ldots, f_s$ in \eqref{eq:sys},
	we have
	\begin{align}\label{eq:f}
		f_i \;=\;(\beta_{i1}-\alpha_{i1}) \left(\kappa_1 \prod\limits_{k=1}^s x_k^{\alpha_{k1}}+\lambda \kappa_2 \prod\limits_{k=1}^s x_k^{\alpha_{k2}} \right), \;\; i=1, \ldots, s.
	\end{align}
	We define the steady-state system augmented with the conservation laws:
	\begin{align}
		h_1 ~&:=~ f_1\;=\;(\beta_{11}-\alpha_{11}) \left(\kappa_1 \prod\limits_{k=1}^s x_k^{\alpha_{k1}}+\lambda \kappa_2 \prod\limits_{k=1}^s x_k^{\alpha_{k2}} \right) ,\label{def_h1} \\
		h_i ~&:= ~(\beta_{i1}-\alpha_{i1}) x_1 - (\beta_{11}-\alpha_{11}) x_i - c_{i-1}, \;\; i=2, \ldots, s.\label{eq:h}
	\end{align}
	We solve $x_i$ from $h_i=0$, and  we get
	\begin{align}\label{eq:ci}
		x_i=\dfrac{(\beta_{i1}-\alpha_{i1}) x_1-c_{i-1}}{\beta_{11}-\alpha_{11}},\; i=2, \ldots, s.
	\end{align}
	We introduce a new variable $z$ and a new parameter $\mu_1$ such that
	\begin{align}\label{eq:1}
		x_1=(\beta_{11}-\alpha_{11}) (z+\mu_1).
	\end{align}
	Then, the conservation laws (i.e., $h_i=0$) can be written as
	\begin{align}\label{eq:2}
		x_i=(\beta_{i1}-\alpha_{i1}) (z+\mu_i),\ \text{where}\;\; \mu_i := \mu_1 - \dfrac{c_{i-1}}{(\beta_{11}-\alpha_{11}) (\beta_{i1}-\alpha_{i1})}. \; 
	\end{align}
	If $h_1=0$, then by \eqref{def_h1} we have $\prod\limits_{k=1}^s x_k^{\alpha_{k1}-\alpha_{k2}}=-\dfrac{\lambda \kappa_2}{\kappa_1} $. So, \begin{align}\label{mihuo1001} \sum\limits_{k=1}^s (\alpha_{k1}-\alpha_{k2})\ln x_k=\ln(-\dfrac{\lambda \kappa_2}{\kappa_1}).\end{align}
	Notice that we can replace $x_i$ with \eqref{eq:2}. So, we define the left hand side of \eqref{mihuo1001} as a new univariate function $g(z)$:
	\begin{align}\label{eq:3.5}
		g(z)\; :&=\; \sum\limits^{s}_{i=1} (\alpha_{i1}-\alpha_{i2}) \ln(\beta_{i1}-\alpha_{i1}) (z+\mu_i).
	\end{align}
	For $i\in S_1,\ S_4$, we define $d_i:=\; \mu_i$. For $i\in S_2,\ S_3$, we define $d_i:=\; -\mu_i$. Recall that we have defined the notions $a_i,\ \gamma_i$ in \eqref{eq:ar}. So, we have 
	\begin{align}\label{eq:4}
		g(z)\; 
		&=\sum\limits_{i\in S_1}a_{i}\ln\gamma_{i}(z+d_{i})-\sum\limits_{i\in S_2}a_{i}\ln\left(-\gamma_{i}(z-d_{i})\right)\notag \\
		&+\sum\limits_{i\in S_3}a_{i}\ln\left(-\gamma_{i}(z-d_{i})\right)-\sum\limits_{i\in S_4}a_{i}\ln\gamma_{i}(z+d_{i}).
	\end{align}
	Notice that the domain of the function $g(z)$ is  $I\; :=\; (\mathcal{L},\; \mathcal{R})$, where
	\begin{align}\label{eq:5}
		\mathcal{L}&:=\begin{cases}
			\max\{-d_{i}\}_{i\in S_1\cup S_4},\  &S_1\cup S_4 \neq \emptyset\   \\
			-\infty, &S_1=S_4=\emptyset
		\end{cases}, \notag \\
		\mathcal{R}&:=\begin{cases}
			\min\{d_{i}\}_{i\in S_2\cup S_3},\ &S_2\cup S_3 \neq \emptyset\   \\
			+\infty, &S_2=S_3=\emptyset
		\end{cases}.
	\end{align}
	By \eqref{eq:3.5} and \eqref{eq:4}, we have
	\begin{align}\label{def dg}
		\dfrac{dg}{dz}(z)=\sum\limits^s_{i=1} \dfrac{\alpha_{i1}-\alpha_{i2}}{ z+\mu_i}
		=\sum\limits_{i\in S_1}\dfrac{a_{i}}{z+d_{i}}+\sum\limits_{i\in S_2}\dfrac{a_{i}}{-z+d_{i}}-\sum\limits_{i\in S_3}\dfrac{a_{i}}{-z+d_{i}}-\sum\limits_{i\in S_4}\dfrac{a_{i}}{z+d_{i}}.
	\end{align}

	\begin{lemma}\label{diff di}
		Given a network $G$ \eqref{eq:2network} with a one-dimensional stoichiometric subspace, suppose $cap_{pos}(G)<+\infty$. Let $g(z)$ and $I$ be the function and the interval  defined as in \eqref{eq:4} and \eqref{eq:5}. Then, $G$ admits multistability iff (if and only if) there exist $\{d_i\}^s_{i=1}\subset \mathbb{R}$ and $K\in\mathbb{R}$ such that the equation $g(z)=K$ has at least $2$ solutions $z_1\ \text{and}\ z_2$ in $I$ satisfying $\dfrac{dg}{dz}(z_1)<0\ \text{and}\ \dfrac{dg}{dz}(z_2)<0$, where these $d_i$'s are distinct from each other.
	\end{lemma}

	\begin{lemma}\label{g' trans}
		Given a network $G$ \eqref{eq:2network} with a one-dimensional stoichiometric subspace, suppose $cap_{pos}(G)<+\infty$. Let $g(z)$ and $I=(\mathcal{L},\; \mathcal{R})$ be defined as in \eqref{eq:4} and \eqref{eq:5}. 
		
		\begin{enumerate}[{(i)}]
			\item\label{lemma4.2_1} For any given $\{d_i\}^s_{i=1}\subset \mathbb{R}$, if $\lim\limits_{z \to \mathcal{L}^{+}}g(z)=+\infty,\; \lim\limits_{z \to \mathcal{R}^{-}}g(z)=-\infty,\; \lim\limits_{z \to \mathcal{L}^{+}}\dfrac{dg}{dz}(z)=-\infty,\;\text{and} \lim\limits_{z \to \mathcal{R}^{-}}\dfrac{dg}{dz}(z)=-\infty$, then there exists $K\in\mathbb{R}$ such that the equation $g(z)=K$ has at least $2$ solutions $z_1\ \text{and}\ z_2$ in $I$ satisfying $\dfrac{dg}{dz}(z_1)<0\ \text{and}\ \dfrac{dg}{dz}(z_2)<0$ iff there exists $\widetilde{z}\in I$ such that $\dfrac{dg}{dz}(\widetilde{z})>0$.

			\item\label{lemma4.2_2} For any given $\{d_i\}^s_{i=1}\subset \mathbb{R}$, if $\lim\limits_{z \to \mathcal{L}^{+}}g(z)=-\infty,\; \lim\limits_{z \to \mathcal{R}^{-}}g(z)=-\infty,\; \lim\limits_{z \to \mathcal{L}^{+}}\dfrac{dg}{dz}(z)=+\infty,\;\text{and} \lim\limits_{z \to \mathcal{R}^{-}}\dfrac{dg}{dz}(z)=-\infty$, then 
			 there exists $K\in\mathbb{R}$ such that the equation $g(z)=K$ has at least $2$ solutions $z_1\ \text{and}\ z_2$ in $I$ satisfying $\dfrac{dg}{dz}(z_1)<0\ \text{and}\ \dfrac{dg}{dz}(z_2)<0$ iff there exist $\widetilde{z}_1,\; \widetilde{z}_2\in I$ $(\widetilde{z}_1<\widetilde{z}_2)$ such that $$\dfrac{dg}{dz}(\widetilde{z}_1)<0,\; \text{and}\;\dfrac{dg}{dz}(\widetilde{z}_2)>0.  $$
			
			
			\item\label{lemma4.2_3} For any given $\{d_i\}^s_{i=1}\subset \mathbb{R}$, if $\lim\limits_{z \to \mathcal{L}^{+}}g(z)=-\infty,\; \lim\limits_{z \to \mathcal{R}^{-}}g(z)=+\infty,\; \lim\limits_{z \to \mathcal{L}^{+}}\dfrac{dg}{dz}(z)=+\infty,\;\text{and} \lim\limits_{z \to \mathcal{R}^{-}}\dfrac{dg}{dz}(z)=+\infty$, and if there exists $K\in\mathbb{R}$ such that the equation $g(z)=K$ has at least $2$ solutions $z_1\ \text{and}\ z_2$ in $I$ satisfying $\dfrac{dg}{dz}(z_1)<0\ \text{and}\ \dfrac{dg}{dz}(z_2)<0$, then 
there exist $\widetilde{z}_1,\; \widetilde{z}_2,\; \widetilde{z}_3\in I$ ($\widetilde{z}_1<\widetilde{z}_2<\widetilde{z}_3$) such that $$\dfrac{dg}{dz}(\widetilde{z}_1)<0,\; \dfrac{dg}{dz}(\widetilde{z}_2)>0,\; \text{and}\;
				\dfrac{dg}{dz}(\widetilde{z}_3)<0.$$
		\end{enumerate}
	\end{lemma}
		
	
	\subsection{Proof of Theorem \ref{thm:main} \ (a)}
	\textbf{$\Leftarrow)$} First, we prove the sufficiency. 
	Without loss of generality, we 
	assume that 
	\begin{align}\label{chongfen 1}
		\sum\limits_{i\in S_1}a_i > \min\limits_{i\in S_4} \{a_i \}.
	\end{align}
	(If $\sum\limits_{i\in S_2}a_i > \min\limits_{i\in S_3} \{a_i \}$, one can similarly prove the network $G$ admits multistability.)
	By Lemma \ref{diff di} and Lemma \ref{g' trans} (i), we only need to find  $\{d_i\}^s_{i=1}\subset \mathbb{R}$ such that $\lim\limits_{z \to \mathcal{L}^{+}}g(z)=+\infty,\; \lim\limits_{z \to \mathcal{R}^{-}}g(z)=-\infty,\; \lim\limits_{z \to \mathcal{L}^{+}}\dfrac{dg}{dz}(z)=-\infty, \lim\limits_{z \to \mathcal{R}^{-}}\dfrac{dg}{dz}(z)=-\infty$, and there exists  $\widetilde{z}\in I$ satisfying $\dfrac{dg}{dz}(\widetilde{z})>0$.
	First, for any $i\in S_1$, we let  $d_i =\sigma_1$. Similarly, we let $d_i = \sigma_2$ for any $i \in S_2$, and we let  $d_i = \sigma_3$ for any $i \in S_3$.  
	Assume that $\min\limits_{i\in S_4} \{a_i \}=a_{i_0}$, where $i_0\in S_4$. And for any $i\in S_4\backslash \{i_0\}$, we also make all $d_i$'s the same, i.e. we let $d_i=\sigma_4 \ \text{for any}\ i \in S_4\backslash \{i_0\}$.
	Then, by \eqref{eq:4} and \eqref{def dg},  we have \begin{align} \label{mihuo 0.4}
		g(z)&=\sum\limits_{i\in S_1}a_i\ln\gamma_i(z+\sigma_1)-\sum\limits_{i\in S_2}a_i\ln\left(-\gamma_i(z-\sigma_2)\right)\notag \\
		&+\sum\limits_{i\in S_3}a_i\ln\left(-\gamma_i(z-\sigma_3)\right)-a_{i_0}\ln\gamma_{i_0}(z+d_{i_0})-\sum\limits_{i\in S_4\backslash \{i_0\} }a_i\ln\gamma_i(z+\sigma_4), \;\text{and}\notag \\
		\dfrac{dg}{dz}(z)&=\dfrac{\sum\limits_{i\in S_1}a_i}{z+\sigma_1}+\dfrac{\sum\limits_{i\in S_2}a_i}{-z+\sigma_2}-\dfrac{\sum\limits_{i\in S_3}a_i}{-z+\sigma_3}-\dfrac{a_{i_0}}{z+d_{i_0}}-\dfrac{\sum\limits_{i\in S_4\backslash\{i_0\}}a_i}{z+\sigma_4}. \end{align}
	So, \begin{align}\label{mihuo 0.5}
		\dfrac{dg}{dz}(0)=\dfrac{\sum\limits_{i\in S_1}a_i}{\sigma_1}+\dfrac{\sum\limits_{i\in S_2}a_i}{\sigma_2}-\dfrac{\sum\limits_{i\in S_3}a_i}{\sigma_3}-\dfrac{a_{i_0}}{d_{i_0}}-\dfrac{\sum\limits_{i\in S_4\backslash\{i_0\}}a_i}{\sigma_4}.
	\end{align}
	Below, we will choose concrete values for $\sigma_1$, $\sigma_2$, $\sigma_3$, $\sigma_4$ and $d_{i_0}$ such that 
	$\dfrac{dg}{dz}(0)>0$.
	We let \begin{align}\label{mihuo 1.5} \sigma_1=\dfrac{\sum\limits_{i\in S_1}a_i+a_{i_0}}{2a_{i_0}},\;\sigma_3=\dfrac{2\sum\limits_{i\in S_3}a_i}{c_0},\sigma_4=\max\left\{\dfrac{2\sum\limits_{i\in S_4\backslash \{i_0\}}a_i}{c_0},\ 2   \right\}\;\text{and} \;d_{i_0}=1,\end{align}
	where \begin{align}\label{mihuo 2} c_0:=\dfrac{\sum\limits_{i\in S_1}a_i}{\sigma_1}-\dfrac{a_{i_0}}{d_{i_0}}=\dfrac{2a_{i_0}\sum\limits_{i\in S_1}a_i}{\sum\limits_{i\in S_1}a_i+a_{i_0}}-a_{i_0}.\end{align}
	Notice that	by (\ref{chongfen 1}), we have $\sum\limits_{i\in S_1}a_i>a_{i_0}$.
	So, by \eqref{mihuo 1.5},  we have $\sigma_1>d_{i_0}$ and by \eqref{mihuo 2}, we have 
	\begin{align}\label{mihuo 3} c_0=\dfrac{a_{i_0}}{\sum\limits_{i\in S_1}a_i+a_{i_0}}(\sum\limits_{i\in S_1}a_i-a_{i_0})>0. \end{align}
	Hence, $\sigma_3>0$ (recall that by \eqref{eq:s} and \eqref{eq:ar}, $a_i>0$ for any $i\not\in S_5$).
	Obviously, we can choose $\sigma_2$ such that $\sigma_2>\sigma_3$.
	Notice that $\dfrac{\sum\limits_{i\in S_2}a_i}{\sigma_2}>0$ (i.e., the second term of $\dfrac{dg}{dz}(0)$ in \eqref{mihuo 0.5} is positive). 
	By \eqref{mihuo 1.5}, we have $\dfrac{c_0}{2}=\dfrac{\sum\limits_{i\in S_3}a_i}{\sigma_3}$, and $\dfrac{c_0}{2}\geqslant\dfrac{\sum\limits_{i\in S_4\backslash \{i_0\}}a_i}{\sigma_4}$. So, by \eqref{mihuo 0.5}, we have 
	\begin{align}\label{mihuo 7} \dfrac{dg}{dz}(0)>\dfrac{\sum\limits_{i\in S_1}a_i}{\sigma_1}-\dfrac{a_{i_0}}{d_{i_0}}-\dfrac{\sum\limits_{i\in S_3}a_i}{\sigma_3}-\dfrac{\sum\limits_{i\in S_4\backslash\{i_0\}}a_i}{\sigma_4}\geqslant c_0-\dfrac{c_0}{2}-\dfrac{c_0}{2}=0.\end{align}
	Notice that for  the interval $I=(\mathcal{L},\; \mathcal{R})$ defined in \eqref{eq:5}, we have  $$\mathcal{L}=\max\{-\sigma_1,\; -d_{i_0}, \ -\sigma_4 \}=-d_{i_0}<0, \;\text{and}\;
	 \mathcal{R}=\min\{\sigma_2,\; \sigma_3\}=\sigma_3>0.$$
 So, we have $0\in I$. 
	By \eqref{mihuo 0.4}, when $z\rightarrow \mathcal{R}^{-}$, we have $\ \dfrac{dg}{dz}(z) \rightarrow -\infty,\ g(z)\rightarrow -\infty$, and when $z\rightarrow \mathcal{L}^{+}$, we have $\ \dfrac{dg}{dz}(z) \rightarrow -\infty,\ g(z)\rightarrow +\infty.$
	Hence, by \eqref{mihuo 7} and by Lemma \ref{diff di} and Lemma \ref{g' trans} (i), the network $G$ admits multistablity.

	\textbf{$\Rightarrow)$} Next, we prove the necessity. Assume that the network $G$ admits multistability. 
 By Lemma \ref{diff di}, $G$ admits multistability iff  there exist $\{d_i\}^s_{i=1}\subset \mathbb{R}$ and $K\in\mathbb{R}$ such that the equation $g(z)=K$ has at least $2$ solutions $z_1\ \text{and}\ z_2$ in the interval $I=(\mathcal{L},\; \mathcal{R})$ defined in \eqref{eq:5} satisfying $\dfrac{dg}{dz}(z_1)<0\ \text{and}\ \dfrac{dg}{dz}(z_2)<0$, where these $d_i$'s are distinct from each other. Below, 
 We prove the conclusion by deducing a contradiction. 
 We assume that 
 \begin{align*}
  \sum\limits_{i\in S_1}a_i \leqslant \min\limits_{i\in S_4} \{a_i \}\ \text{and} \ \sum\limits_{i\in S_2}a_i \leqslant \min\limits_{i\in S_3} \{a_i \} .  
 \end{align*}
 Equivalently, we have
	\begin{align}
		\sum\limits_{i\in S_1}a_i &\leqslant a_{j}, \; \text{for any}\; j\in S_4\; \text{and}\label{jiashe1} \\ \sum\limits_{i\in S_2}a_i &\leqslant a_{j}, \; \text{for any},\; j\in S_3.\label{jiashe2}
	\end{align}
	 By \eqref{eq:5}, there exists $m\in S_1\cup S_4$ such that $\mathcal{L}=-d_m$ and there exists $n\in S_2\cup S_3$ such that $\mathcal{R}=d_n$. Below, we deduce the contradiction  for  four different cases.
	
	\textbf{(Case 1)} Assume that there exist $m\in S_1$ and $n\in S_3$ such that \begin{align}\label{mihuo 10} \mathcal{L}=-d_m,\ \text{and} \ \mathcal{R}=d_n. \end{align}
	By \eqref{eq:4} and \eqref{def dg}, we have $$\lim\limits_{z \to -d^{+}_m}g(z)=-\infty,\ \lim\limits_{z \to d_n^{-}}g(z)=-\infty,\ \lim\limits_{z \to -d_m^{+}}\dfrac{dg}{dz}(z)=+\infty,\;\text{and}\; \lim\limits_{z \to d_n^{-}}\dfrac{dg}{dz}(z)=-\infty.$$
	So, by Lemma \ref{diff di} and Lemma \ref{g' trans} (ii) (2), there exists $z_0\in I=(-d_m,\; d_n)$ such that
	\begin{align}
		\dfrac{dg}{dz}(z_0)&=\sum\limits_{i\in S_1}\dfrac{a_i}{z_0+d_i}+\sum\limits_{i\in S_2}\dfrac{a_i}{-z_0+d_i}-\sum\limits_{i\in S_3}\dfrac{a_i}{-z_0+d_i}-\sum\limits_{i\in S_4}\dfrac{a_i}{z_0+d_i}=0, \;\text{and}\label{a dgz0}\\
		\dfrac{d^2 g}{dz^2}(z_0)&=-\sum\limits_{i\in S_1}\dfrac{a_i}{(z_0+d_i)^2}+\sum\limits_{i\in S_2}\dfrac{a_i}{(-z_0+d_i)^2}-\sum\limits_{i\in S_3}\dfrac{a_i}{(-z_0+d_i)^2}+\sum\limits_{i\in S_4}\dfrac{a_i}{(z_0+d_i)^2}\geqslant 0.\label{a ddgz0}
	\end{align}
 Below, we show that 
 if \eqref{a dgz0} and \eqref{a ddgz0} hold Simultaneously, then there will be a contradiction. 
 By \eqref{eq:5} and \eqref{mihuo 10}, we have
	\begin{align}\label{mihuo 11} \mathcal{R}=d_n<d_i \ \text{for any}\ i \in S_2. \end{align} 
	So, \begin{align*} \dfrac{a_i}{-z_0+d_i}<\dfrac{a_i}{-z_0+d_n}  \ \text{for any}\ i \in S_2. \end{align*} 
	Thus, for the second and the third terms in \eqref{a dgz0}, we have
	\begin{align}\label{mihuo 59} \sum\limits_{i\in S_2}\dfrac{a_i}{-z_0+d_i}-\sum\limits_{i\in S_3}\dfrac{a_i}{-z_0+d_i}<\dfrac{\sum\limits_{i\in S_2}a_i}{-z_0+d_n}- \sum\limits_{i\in S_3}\dfrac{a_i}{-z_0+d_i}.\end{align}
	By the fact that $n\in S_3$ and by \eqref{jiashe2}, we have
	\begin{align}\label{mihuo 60} \dfrac{\sum\limits_{i\in S_2}a_i}{-z_0+d_n}-\sum\limits_{i\in S_3}\dfrac{a_i}{-z_0+d_i}\leqslant \dfrac{\sum\limits_{i\in S_2}a_i}{-z_0+d_n}-\dfrac{a_n}{-z_0+d_n}
 \leqslant 0.\end{align}
So,	by \eqref{mihuo 59} and \eqref{mihuo 60}, we have \begin{align}\label{mihuo 13}
		\sum\limits_{i\in S_2}\dfrac{a_i}{-z_0+d_i}<\sum\limits_{i\in S_3}\dfrac{a_i}{-z_0+d_i}.
	\end{align}
	Then, by (\ref{a dgz0}), for the first and the last terms in \eqref{a dgz0}, we have
	\begin{align}\label{a bds1}
		\sum\limits_{i\in S_1}\dfrac{a_i}{z_0+d_i} > \sum\limits_{i\in S_4}\dfrac{a_i}{z_0+d_i}.
	\end{align}
	Similarly, by \eqref{mihuo 11} and \eqref{jiashe2},   we have
	$$\sum\limits_{i\in S_2}\dfrac{a_i}{(-z_0+d_i)^2}-\sum\limits_{i\in S_3}\dfrac{a_i}{(-z_0+d_i)^2}<\dfrac{\sum\limits_{i\in S_2}a_i}{(-z_0+d_n)^2}- \sum\limits_{i\in S_3}\dfrac{a_i}{(-z_0+d_i)^2}$$ 
	$$\leqslant \dfrac{\sum\limits_{i\in S_2}a_i}{(-z_0+d_n)^2}-\dfrac{a_n}{(-z_0+d_n)^2} \leqslant 0.$$
	Then, by (\ref{a ddgz0}), we have
	\begin{align}\label{a bds2}
		\sum\limits_{i\in S_4}\dfrac{a_i}{(z_0+d_i)^2}> \sum\limits_{i\in S_1}\dfrac{a_i}{(z_0+d_i)^2}.
	\end{align}
	By (\ref{jiashe1}), we have 
	\begin{align}\label{mihuo 63}
		\sum\limits_{i\in S_4}\dfrac{a_i^2}{(z_0+d_i)^2}\geqslant 
  \sum\limits_{i\in S_4}\dfrac{a_i}{(z_0+d_i)^2}\cdot \sum\limits_{i\in S_1}a_i.
	\end{align}
	By (\ref{a bds2}), by using the Cauchy-Schwarz inequality and by \eqref{a bds1}, we have
	$$\sum\limits_{i\in S_4}\dfrac{a_i}{(z_0+d_i)^2}\cdot \sum\limits_{i\in S_1}a_i> \sum\limits_{i\in S_1}\dfrac{a_i}{(z_0+d_i)^2} \cdot \sum\limits_{i\in S_1}a_i \geqslant \left(\sum\limits_{i\in S_1}\dfrac{a_i}{z_0+d_i}\right)^2 >\left(\sum\limits_{i\in S_4}\dfrac{a_i}{z_0+d_i}\right)^2. $$
	So, by \eqref{mihuo 63}, we have $$\sum\limits_{i\in S_4}\dfrac{a_i^2}{(z_0+d_i)^2}> \left(\sum\limits_{i\in S_4}\dfrac{a_i}{z_0+d_i}\right)^2,$$
	which is impossible since
 $\dfrac{a_i}{z_0+d_i}>0$.
	
	\textbf{(Case 2)} Assume that there exist $m\in S_4$ and $n\in S_2$ such that \begin{align}\label{mihuo 21} \mathcal{L}=-d_m,\;\text{and}\;\mathcal{R}=d_n. \end{align}
	Since Case 2 is symmetric with respect to  Case 1, the proof is similar to the proof of Case 1. We omit the details. 
	
	\textbf{(Case 3)} Assume that there exist $m\in S_4$ and $n\in S_3$ such that \begin{align}\label{mihuo 21} \mathcal{L}=-d_m,\;\text{and}\;\mathcal{R}=d_n. \end{align}
	Then, by \eqref{eq:5}, we have $$-d_m=\max\{-d_i\}_{i\in S_1\cup S_4},\;\text{and}\; d_n=\min\{d_i\}_{i\in S_2\cup S_3}.$$
	So, for any $z\in (-d_m,\; d_n)$, we have  \begin{align}\label{mihuo 101}
		\dfrac{a_i}{z+d_i}<\dfrac{a_i}{z+d_m} \ \text{for any}\ i \in S_1,\;\text{and}\;\notag \\
		\dfrac{a_i}{-z+d_i}<\dfrac{a_i}{-z+d_n} \ \text{for any}\ i \in S_2.
	\end{align}
	Since $m\in S_4\;\text{and}\; n\in S_3$, by (\ref{jiashe1})--(\ref{jiashe2}), we have 
	\begin{align}\label{eq:henmihuo}
		a_m \geqslant \sum\limits_{i\in S_1}a_i,\;\text{and}\; a_n\geqslant \sum\limits_{i\in S_2}a_i.
	\end{align}
	Then, by \eqref{def dg}, \eqref{mihuo 101} and \eqref{eq:henmihuo}, we have
	$$\dfrac{dg}{dz}(z)=\sum\limits_{i\in S_1}\dfrac{a_i}{z+d_i}+\sum\limits_{i\in S_2}\dfrac{a_i}{-z+d_i}-\sum\limits_{i\in S_3}\dfrac{a_i}{-z+d_i}-\sum\limits_{i\in S_4}\dfrac{a_i}{z+d_i}$$ $$\leqslant 
	\dfrac{\sum\limits_{i\in S_1}a_i}{z+d_m}+\dfrac{\sum\limits_{i\in S_2}a_i}{-z+d_n}-\dfrac{a_n}{-z+d_n}-\dfrac{a_m}{z+d_m}\leqslant 0.
	$$
	So, $g(z)$ is decreasing in $I$. Thus, $g(z)=0$ has at most one real solution in $I$. On the other hand, by Lemma \ref{diff di}, $g(z)=0$ has at least two real solutions in $I$, which is a contradiction.
	
	\textbf{(Case 4)} Assume that there exist $m\in S_1$ and $n\in S_2$ such that \begin{align}\label{mihuo 21} \mathcal{L}=-d_m,\;\text{and}\; \mathcal{R}=d_n. \end{align}
	By \eqref{eq:4} and \eqref{def dg}, we have $$\lim\limits_{z \to -d^{+}_m}g(z)=-\infty,\ \lim\limits_{z \to d_n^{-}}g(z)=+\infty,\ \lim\limits_{z \to -d_m^{+}}\dfrac{dg}{dz}(z)=+\infty,\;\text{and}\; \lim\limits_{z \to d_n^{-}}\dfrac{dg}{dz}(z)=+\infty.$$
	So, by Lemma \ref{diff di} and Lemma \ref{g' trans} (iii) (2), if $G$ admits multistability, then  there exists $z_0\in I=(-d_m,\; d_n)$ such that
	$$\dfrac{dg}{dz}(z_0)>0,\ \dfrac{d^2 g}{dz^2}(z_0)=0,\;\text{and}\;  \dfrac{d^3 g}{dz^3}(z_0)\leqslant 0.$$
	So, by \eqref{def dg}, we have
	\begin{align}
		\sum\limits_{i\in S_1}\dfrac{a_i}{z_0+d_i}+\sum\limits_{i\in S_2}\dfrac{a_i}{-z_0+d_i}&>\sum\limits_{i\in S_3}\dfrac{a_i}{-z_0+d_i}+\sum\limits_{i\in S_4}\dfrac{a_i}{z_0+d_i}, \label{d dgz0} \\
		\sum\limits_{i\in S_1}\dfrac{a_i}{(z_0+d_i)^2}+\sum\limits_{i\in S_3}\dfrac{a_i}{(-z_0+d_i)^2}&=\sum\limits_{i\in S_2}\dfrac{a_i}{(-z_0+d_i)^2}+\sum\limits_{i\in S_4}\dfrac{a_i}{(z_0+d_i)^2}, \;\text{and}\label{d ddgz0} \\
		\sum\limits_{i\in S_1}\dfrac{a_i}{(z_0+d_i)^3}+\sum\limits_{i\in S_2}\dfrac{a_i}{(-z_0+d_i)^3}&\leqslant \sum\limits_{i\in S_3}\dfrac{a_i}{(-z_0+d_i)^3}+\sum\limits_{i\in S_4}\dfrac{a_i}{(z_0+d_i)^3}. \label{d dddgz0} 
	\end{align}
	Note that by (\ref{d dddgz0}), one of the following two equations must hold
	\begin{align} \sum\limits_{i\in S_1}\dfrac{a_i}{(z_0+d_i)^3}&\leqslant \sum\limits_{i\in S_4}\dfrac{a_i}{(z_0+d_i)^3},\label{mihuo 50} \\  \sum\limits_{i\in S_2}\dfrac{a_i}{(-z_0+d_i)^3}&\leqslant \sum\limits_{i\in S_3}\dfrac{a_i}{(-z_0+d_i)^3}.\label{mihuo 51} \end{align}
	Without loss of generality, we assume \eqref{mihuo 50} holds (if \eqref{mihuo 51} holds, we can prove the conclusion similarly).  Below, we prove that \eqref{d dgz0}, \eqref{d ddgz0} and \eqref{mihuo 50} can not hold simultaneously  by the following four steps.
	
	\textbf{(Step 1)} In this step, we prove that \eqref{d dgz0} and \eqref{d ddgz0} imply 
	\begin{align}\label{d **}
		\sum\limits_{i\in S_1}\dfrac{a_i}{(z_0+d_i)^2}> \sum\limits_{i\in S_4}\dfrac{a_i}{(z_0+d_i)^2}.
	\end{align}
	We prove the conclusion  by deducing a contradiction. Assume that \begin{align} \label{mihuo 30} \sum\limits_{i\in S_1}\dfrac{a_i}{(z_0+d_i)^2}\leqslant \sum\limits_{i\in S_4}\dfrac{a_i}{(z_0+d_i)^2}. \end{align}
	Then, by Cauchy inequality and by \eqref{mihuo 30}, we have
	\begin{align}\label{mihuo 30.1} \left(\sum\limits_{i\in S_1}\dfrac{a_i}{z_0+d_i}\right)^2\leqslant \sum\limits_{i\in S_1}\dfrac{a_i}{(z_0+d_i)^2}\cdot \sum\limits_{i\in S_1}a_i\leqslant \sum\limits_{i\in S_4}\dfrac{a_i}{(z_0+d_i)^2}\cdot \sum\limits_{i\in S_1}a_i.\end{align}
	By (\ref{jiashe1}), we have 
	\begin{align}\label{mihuo 30.2} \sum\limits_{i\in S_4}\dfrac{a_i}{(z_0+d_i)^2}\cdot \sum\limits_{i\in S_1}a_i\leqslant \sum\limits_{i\in S_4}\dfrac{a_i^2}{(z_0+d_i)^2}\leqslant \left(\sum\limits_{i\in S_4}\dfrac{a_i}{z_0+d_i}\right)^2. \end{align}
	Notice that $\dfrac{a_i}{z_0+d_i}>0$ for any $z_0\in I$. So, by \eqref{mihuo 30.1} and \eqref{mihuo 30.2}, we have
	$$\sum\limits_{i\in S_1}\dfrac{a_i}{z_0+d_i}\leqslant \sum\limits_{i\in S_4}\dfrac{a_i}{z_0+d_i}.$$
	Hence, by (\ref{d dgz0}), we have
	\begin{align}\label{d lue1}
		\sum\limits_{i\in S_2}\dfrac{a_i}{-z_0+d_i}> \sum\limits_{i\in S_3}\dfrac{a_i}{-z_0+d_i}.
	\end{align}
	On the other hand, by (\ref{d ddgz0}) and \eqref{mihuo 30}, we have
	\begin{align}\label{eq:add_1shang}
	\sum\limits_{i\in S_2}\dfrac{a_i}{(-z_0+d_i)^2}\leqslant \sum\limits_{i\in S_3}\dfrac{a_i}{(-z_0+d_i)^2}
	.
 \end{align}
	So, by Cauchy inequality and by \eqref{eq:add_1shang}, we have \begin{align}\label{add_1}
 \left(\sum\limits_{i\in S_2}\dfrac{a_i}{-z_0+d_i}\right)^2\leqslant \sum\limits_{i\in S_2}\dfrac{a_i}{(-z_0+d_i)^2}\cdot \sum\limits_{i\in S_2}a_i\leqslant \sum\limits_{i\in S_3}\dfrac{a_i}{(-z_0+d_i)^2}\cdot \sum\limits_{i\in S_2}a_i.\end{align}
	By (\ref{jiashe2}) and by the fact $a_{i_2}\leqslant a_i$ for any $i\in S_3$, we have
	\begin{align}\label{add_2}
 \sum\limits_{i\in S_3}\dfrac{a_i}{(-z_0+d_i)^2}\cdot \sum\limits_{i\in S_2}a_i \leqslant \sum\limits_{i\in S_3}\dfrac{a_i^2}{(-z_0+d_i)^2} \leqslant \left(\sum\limits_{i\in S_3}\dfrac{a_i}{-z_0+d_i}\right)^2.\end{align}
	So, by \eqref{add_1} and \eqref{add_2}, we have 
	$$\sum\limits_{i\in S_2}\dfrac{a_i}{-z_0+d_i}\leqslant \sum\limits_{i\in S_3}\dfrac{a_i}{-z_0+d_i}.$$
	This is a contradiction to (\ref{d lue1}). Therefore, the inequality \eqref{d **} must hold.
	
	\textbf{(Step 2)} In this step, we show that by (\ref{mihuo 50}) and (\ref{d **}), we can construct $d_0,\ a_0\in \mathbb{R}$ such that 	\begin{align}
		\sum\limits_{i\in S_1}\dfrac{a_i}{(z_0+d_i)^2}> \dfrac{a_0}{(z_0+d_0)^2},\;\text{and}\label{quan 1} \\
		\sum\limits_{i\in S_1}\dfrac{a_i}{(z_0+d_i)^3}\leqslant \dfrac{a_0}{(z_0+d_0)^3}.\label{quan 2}
	\end{align}
	Let \begin{align}	d_0\; :=\; \left(\sum\limits_{i\in S_4}\dfrac{a_i}{(z_0+d_i)^2}\right)\bigg/  \left(\sum\limits_{i\in S_4}\dfrac{a_i}{(z_0+d_i)^3}\right)-z_0,\label{mihuo 32} \\
		a_0\; :=\; \left(\sum\limits_{i\in S_4}\dfrac{a_i}{(z_0+d_i)^2}\right)^3 \bigg/  \left(\sum\limits_{i\in S_4}\dfrac{a_i}{(z_0+d_i)^3}\right)^2. \label{mihuo 33}
	\end{align}
	It is straightforward to check that
	\begin{align}\label{d0 bubian}
		\sum\limits_{i\in S_4}\dfrac{a_i}{(z_0+d_i)^2}=\dfrac{a_0}{(z_0+d_0)^2},\;\text{and}\; \sum\limits_{i\in S_4}\dfrac{a_i}{(z_0+d_i)^3}=\dfrac{a_0}{(z_0+d_0)^3}.
	\end{align}
	So, by (\ref{mihuo 50}) and (\ref{d **}), we have \eqref{quan 1} and \eqref{quan 2} hold.
	
	\textbf{(Step 3)} In this step, we prove that
	\begin{align}\label{a0 is max}
		a_0\geqslant \sum\limits_{i\in S_1}a_i.
	\end{align}
	By \eqref{mihuo 33}, we only need to show that  \begin{align}\label{mihuo 53} \left(\sum\limits_{i\in S_4}\dfrac{a_i}{(z_0+d_i)^2}\right)^3\geqslant \left(\sum\limits_{i\in S_4}\dfrac{a_i}{(z_0+d_i)^3}\right)^2\cdot \sum\limits_{i\in S_1}a_i. \end{align}
	Notice that by Cauchy inequality, we have \begin{align}\label{mihuo 54} \left(\sum\limits_{i\in S_4}\dfrac{a_i}{(z_0+d_i)^2}\right)^3\cdot &\sum\limits_{i\in S_4}\dfrac{a_i}{(z_0+d_i)^4}= \left(\sum\limits_{i\in S_4}\dfrac{a_i}{(z_0+d_i)^2}\right)^2\cdot \sum\limits_{i\in S_4}\dfrac{a_i}{(z_0+d_i)^2} \cdot \sum\limits_{i\in S_4}\dfrac{a_i}{(z_0+d_i)^4}\\ \notag
 &\geqslant \left(\sum\limits_{i\in S_4}\dfrac{a_i}{(z_0+d_i)^2}\right)^2\cdot \left(\sum\limits_{i\in S_4}\dfrac{a_i}{(z_0+d_i)^3}\right)^2.
 \end{align} 
	By \eqref{jiashe1}, we have \begin{align}\label{mihuo 55} \left(\sum\limits_{i\in S_4}\dfrac{a_i}{(z_0+d_i)^2}\right)^2\geqslant \sum\limits_{i\in S_4}\dfrac{a_i^2}{(z_0+d_i)^4} \geqslant \sum\limits_{i\in S_4}\dfrac{a_i}{(z_0+d_i)^4}\cdot \sum\limits_{i\in S_1}a_i.\end{align}
	So, by \eqref{mihuo 54} and \eqref{mihuo 55}, one can see that  \eqref{mihuo 53} holds. 

	\textbf{(Step 4)} In this step, we show that \eqref{quan 1} and \eqref{quan 2} imply  two contradictory inequalities \eqref{maodun1} and \eqref{mihuo 57}.
	By Cauchy inequality and by (\ref{quan 2}), we have
	\begin{align}\label{maodun0}
		\left(\sum\limits_{i\in S_1}\dfrac{a_i}{(z_0+d_i)^2}\right)^2\leqslant \sum\limits_{i\in S_1}\dfrac{a_i}{(z_0+d_i)^3}\cdot \sum\limits_{i\in S_1}\dfrac{a_i}{z_0+d_i}\leqslant \dfrac{a_0}{(z_0+d_0)^3} \cdot \sum\limits_{i\in S_1}\dfrac{a_i}{z_0+d_i}.
	\end{align}
	We multiply the left-hand side and the right-hand side of \eqref{maodun0} by $\sum\limits_{i\in S_1}a_i$, and we get \begin{align}\label{mihuo 37}
		\left(\sum\limits_{i\in S_1}\dfrac{a_i}{(z_0+d_i)^2}\right)^2 \cdot \sum\limits_{i\in S_1} a_i \leqslant \dfrac{a_0}{(z_0+d_0)^3} \cdot \sum\limits_{i\in S_1}\dfrac{a_i}{z_0+d_i} \cdot \sum\limits_{i\in S_1}a_i.
	\end{align}
	Note that by Cauchy inequality, we have \begin{align}\label{mihuo 38}
		\sum\limits_{i\in S_1}\dfrac{a_i}{(z_0+d_i)^2}\cdot \left(\sum\limits_{i\in S_1}\dfrac{a_i}{z_0+d_i}\right)^2 \leqslant \left(\sum\limits_{i\in S_1}\dfrac{a_i}{(z_0+d_i)^2}\right)^2 \cdot \sum\limits_{i\in S_1} a_i . 
	\end{align}
	By (\ref{a0 is max}), we have $\sum\limits_{i\in S_1}a_i \leqslant a_0 $. So, 
	\begin{align}\label{mihuo 39}
		\dfrac{a_0}{(z_0+d_0)^3} \cdot \sum\limits_{i\in S_1}\dfrac{a_i}{z_0+d_i}\cdot \sum\limits_{i\in S_1}a_i\leqslant \dfrac{a_0^2}{(z_0+d_0)^3}\cdot \sum\limits_{i\in S_1}\dfrac{a_i}{z_0+d_i}.
	\end{align}
	Thus, by \eqref{mihuo 37}--\eqref{mihuo 39}, we have
	$$\sum\limits_{i\in S_1}\dfrac{a_i}{(z_0+d_i)^2}\cdot \left(\sum\limits_{i\in S_1}\dfrac{a_i}{z_0+d_i}\right)^2 \leqslant \dfrac{a_0^2}{(z_0+d_0)^3}\cdot \sum\limits_{i\in S_1}\dfrac{a_i}{z_0+d_i}.$$ Note that $\dfrac{a_i}{z_0+d_i}>0\ (i\in S_1)$ for any $z_0\in I$. So, we have
	$$\sum\limits_{i\in S_1}\dfrac{a_i}{(z_0+d_i)^2}\cdot \sum\limits_{i\in S_1}\dfrac{a_i}{z_0+d_i}\leqslant \dfrac{a_0^2}{(z_0+d_0)^3}.$$
	Then, by (\ref{quan 1}), we have
	\begin{align}\label{maodun1}
		\sum\limits_{i\in S_1}\dfrac{a_i}{z_0+d_i}< \dfrac{a_0}{z_0+d_0}.
	\end{align}
	On the other hand, by (\ref{quan 1}) and (\ref{maodun0}), we have
	$$\dfrac{a_0^2}{(z_0+d_0)^4}<\left(\sum\limits_{i\in S_1}\dfrac{a_i}{(z_0+d_i)^2}\right)^2 \leqslant \dfrac{a_0}{(z_0+d_0)^3}\cdot \sum\limits_{i\in S_1}\dfrac{a_i}{z_0+d_i}.$$	
	Note that by   \eqref{mihuo 32} and \eqref{mihuo 33}, we have  $z_0+d_0>0$ and $a_0>0$. So, 
	\begin{align}\label{mihuo 57} \dfrac{a_0}{z_0+d_0}< \sum\limits_{i\in S_1}\dfrac{a_i}{z_0+d_i},\end{align}
	which is a contradiction to (\ref{maodun1}). So far, we have deduced the contradiction for the last case, and we complete the proof.
	
	\subsection{Proof of Theorem \ref{thm:main}\;(b)}
	In this case, there are exactly three of the four sets $S_1,\; S_2,\; S_3,\; S_4$ are non-empty. Below, we successively prove Theorem \ref{thm:main} (b) (1)-(4).
	\subsubsection{Proof of Theorem \ref{thm:main} (b) (1)}\label{134}
	According to the hypothesis of Theorem \ref{thm:main} (b) (1), we assume that $S_1$, $S_3$ and $S_4$ are non-empty.
	By \eqref{eq:4}, we have
	\begin{align}
		\label{g134} g(z)\;&= \;  \sum\limits_{i \in S_1}a_i\ln(z+d_i)
		+\sum\limits_{i \in S_3}a_i\ln(-z+d_i)
		-\sum\limits_{i \in S_4}a_i\ln(z+d_i),\\
		\label{dg134} \dfrac{dg}{dz}(z)\; &= \; \sum\limits_{i \in S_1}\dfrac{a_i}{z+d_i}
		-\sum\limits_{i \in S_3}\dfrac{a_i}{-z+d_i}
		-\sum\limits_{i \in S_4}\dfrac{a_i}{z+d_i}.
	\end{align}
	
	$\Rightarrow)$ First, we prove the sufficiency.
	By \cite[Theorem 3.4]{txzs2020}, if the network G admits multistability, then $cap_{pos}(G)\geq 3$. By \cite[Theorem 6.1 (c)]{linkexin}, if $cap_{pos}(G)\geq 3$, then we have $\sum\limits_{i\in S_1}a_i>\min\limits_{i\in S_{4}}\{a_i\}$. 
	
	$\Leftarrow)$ Next, we prove the necessity. Assume that
	\begin{align}\label{eq:468}
		\sum\limits_{i\in S_1}a_i>\min\limits_{i\in S_{4}}\{a_i\}.
	\end{align}    
 The goal is to prove that $G$ admits multistability. Assume that $a_{p} = \min\limits_{i \in S_4} \{a_i\}$, where $p\in S_4$. 
	First, we let $d_{p} = 0$. 
 Then, for any $i\in S_1$, we let $d_i = d$, and for any $i \in S_3$, we let $d_i = 1$. For any $S_4\setminus{\{p\}}$, we also make all $d_i$'s the same, i.e. we let $d_i = e$ for any $i \in S_4\setminus{\{p\}}$. Notice that  $d$ and $e$ are two positive parameters, and we will choose proper values for them later.
 By \eqref{dg134}, we have 
 \begin{align}\label{eq:dg135}
 \dfrac{dg}{dz}(z) = \dfrac{\sum\limits_{i \in S_1}a_i}{z+d}-\dfrac{\sum\limits_{i \in S_3}a_i}{-z+1}-\dfrac{a_{p}}{z}-\dfrac{\sum\limits_{i \in S_4\setminus{\{p\}}}a_i}{z+e}.
 \end{align}
	Note that the interval  $I$ defined in \eqref{eq:5} is $(0,\;1)$. By \eqref{g134} and \eqref{dg134}, we have $\lim\limits_{z \to 0^{+}}g(z)=+\infty,\; \lim\limits_{z \to 1^{-}}g(z)=-\infty,\; \lim\limits_{z \to 0^{+}}\dfrac{dg}{dz}(z)=-\infty,\;\text{and} \lim\limits_{z \to 1^{-}}\dfrac{dg}{dz}(z)=-\infty$. By Lemma \ref{diff di} and Lemma 4.2 (\ref{lemma4.2_1}), we only need to choose proper positive numbers $d$ and $e$ such that there exists $\widetilde{z}\in I$ satisfying  $\dfrac{dg}{dz}(\widetilde{z})>0$.
	Below, we complete the proof by the following  two steps. 
 
	(Step 1)Let \begin{align}\label{hz134}
		h(z) := \dfrac{\sum\limits_{i \in S_1}a_i}{z+d}-\dfrac{\sum\limits_{i \in S_3}a_i}{-z+1}-\dfrac{a_{p}}{z}.
	\end{align}
	Notice that by \eqref{eq:dg135}, $h(z)=\dfrac{dg}{dz}(z)+\dfrac{\sum\limits_{i \in S_4\setminus{\{p\}}}a_i}{z+e}$. In this step, we show that we can choose $d>0$ such that there exists $\widetilde{z} \in I$ satisfying $h(\widetilde{z}) > 0$.
	We solve $d$ from $h(z)>0$ by \eqref{hz134}, and we get
	\begin{align}\label{d134}
		d < \dfrac{{\mathcal N}(z)}{{\mathcal D}(z)},
	\end{align} 
 where  ${\mathcal D}(z):=\dfrac{\sum\limits_{i \in S_3}a_i}{-z+1}+\dfrac{a_{p}}{z}$ and ${\mathcal N}(z):=\sum\limits_{i \in S_1}a_i-\dfrac{z}{-z+1}\sum\limits_{i \in S_3}a_i-a_{p}$.
	Notice that by \eqref{eq:468}, we have $\sum\limits_{i \in S_1}a_i>a_{p}$, and so,  
	\begin{align}\label{num>0}
		\lim\limits_{z \to 0}{\mathcal N}(z)>0.
	\end{align}
	So, there exists $\widetilde{z} \in I=(0,1)$, such that the ${\mathcal N}(\widetilde{z})>0$.
	Notice that for any $z\in I =(0,\;1)$, we have 
 ${\mathcal D}(z)>0$. 
	Therefore, $\dfrac{{\mathcal N}(\widetilde{z})}{{\mathcal D}(\widetilde{z})}>0$.
	Then, we can choose an appropriate positive number $d$ such that  $d<\dfrac{{\mathcal N}(\widetilde{z})}{{\mathcal D}(\widetilde{z})}$, i.e., $h(\widetilde{z})>0$.
	
	(Step 2) In this step, we prove that we can choose $e>0$ such that  $\dfrac{dg}{dz}(\widetilde{z}) > 0$.
	In fact, let $$e = 2~\dfrac{\sum\limits_{i \in S_4\setminus{\{p\}}}a_i}{h(\widetilde{z})}.$$
	Then, $$\dfrac{dg}{dz}(\widetilde{z}) = h(\widetilde{z})-\dfrac{\sum\limits_{i \in S_4\setminus{\{p\}}}a_i}{\widetilde{z}+e}>h(\widetilde{z})-\dfrac{\sum\limits_{i \in S_4\setminus{\{p\}}}a_i}{e}=h(\widetilde{z})-\dfrac{h(\widetilde{z})}{2}>0.$$
	\subsubsection{Proof of Theorem \ref{thm:main} (b) (2):\; $S_2$, $S_3$ and $S_4$ are non-empty}\label{234}
 According to the hypothesis of Theorem \ref{thm:main} (b) (2), we assume that $S_2$, $S_3$ and $S_4$ are non-empty.
	By \eqref{eq:4}, we have
	\begin{align*}
		g(z)\;&= \;  
		-\sum\limits_{i \in S_2}a_i\ln(-z+d_i)
		+\sum\limits_{i \in S_3}a_i\ln(-z+d_i)
		-\sum\limits_{i \in S_4}a_i\ln(z+d_i).
	\end{align*}
	Define \begin{align*}
		\widetilde{g}(z)\;:=\; -g(-z)\;=\;  
		\sum\limits_{i \in S_2}a_i\ln(z+d_i)
		-\sum\limits_{i \in S_3}a_i\ln(z+d_i)
		+\sum\limits_{i \in S_4}a_i\ln(-z+d_i).
	\end{align*}
	 Notice that $\dfrac{d\widetilde{g}}{dz}(z) = \dfrac{dg}{dz}(-z)$. Let $I^* = \{-z|z\in I\}$. Then, there exist $z_1,\;z_2\in I$ such that  $g(z_i)=0,\;\text{and}\; \dfrac{dg}{dz}(z_i)<0$ ($i=1,2$) if and only if there exist $z_1^*,\;z_2^*\in I^*$ such that  $\widetilde{g}(z_i^*)=0,\; \text{and}\; \dfrac{d\widetilde{g}}{dz}(z_i^*)<0$ ($i=1,2$). Note that by the proof of Theorem \ref{thm:main} (b) (1), there exist $z_1^*,\;z_2^*\in I^*$ such that  $\widetilde{g}(z_i^*)=0,\;  \text{and } \dfrac{d\widetilde{g}}{dz}(z_i^*)<0$ ($i=1, 2$) if and only if $
	\sum\limits_{i\in S_2}a_{i}>\min\limits_{i\in S_{3}}\{a_i\}.$ So, by Lemma \ref{diff di}, $G$ admits multistability if and only if $
	\sum\limits_{i\in S_2}a_{i}>\min\limits_{i\in S_{3}}\{a_i\}.$
	\subsubsection{Proof of Theorem \ref{thm:main} (b) (3)}\label{123}

	First, we give some lemmas (Lemma \ref{simplify1}--Lemma \ref{e2<}). Since the proofs of these lemmas are elementary, we put them in the supplementary materials \footnote{https://github.com/65536-1024/one-dim}. 
	\begin{lemma}\label{simplify1}
		For any $\beta_1,\;\beta_2,\;e_1,\;e_2,\;x_1,\;x_2,\;x_3\in R$ satisfying
\begin{align}
\label{7number_1}	&\beta_1>0,\; \beta_2>0,\;\\ 
\label{7number_2} &x_i+e_j>0 \;(i=1,\;2,\;3,\;j=1,\;2), \;\text{and} \\
\label{7number_3} &e_1\neq e_2,
		\end{align} where exist $a,\;b,\;c\in R$ such that 
		\begin{align}
			\label{eq:fourzerolemma'}\dfrac{\beta_1}{x_i+e_1}+\dfrac{\beta_2}{x_i+e_2}  &= a + \dfrac{c}{x_i+b}\;(i=1,\;2,\;3),\end{align}
		where\begin{align}\label{abc}\;\;\;\;
			a>0,\; b>\min\{e_1,\;e_2\},\;\text{and}\; \min\{\beta_1, \;\beta_2\}<c<\beta_1+\beta_2.
		\end{align}
	\end{lemma}
	\begin{lemma}\label{simplify3}Let $G(z) := \sum\limits_{i=1}\limits^{n}\dfrac{a_{i}}{z+d_{i}}$, where $d_i\in R$, and $a_i>0$. Let $M := \min\limits_{i\in\{1,\;\cdots,\;n\}}\{d_i\}$. Then, for any three different numbers $z_1, z_2, z_3$ satisfying $z_j > -M\;(j=1,\;2,\;3)$, there exist $A,\;D,\;\theta\in R$ such that
		\begin{align}\label{elimi1}
			G(z_j) = \dfrac{A}{z_j+D}+\theta \;( j = 1,\;2,\;3),
		\end{align}
		where 
		\begin{align}
			\min\limits_{i\in\{1,\;\cdots,\;n\}}\{a_i\}\leq A\leq\sum\limits_{i=1}\limits^{n}a_i,\;D\geq M,\;\text{and}\;\theta\geq0
		\end{align}
	\end{lemma}

	\begin{lemma}\label{finalineq}
		For any two sequences $\{a_i\}_{i=0}^{n}$, $\{e_i\}_{i=0}^{n}$ satisfying 
		\begin{align}
			\label{lemfinalcon} a_i>1(i=1,\;\cdots,\;n),\;e_i>1(i=1,\;\cdots,\;n),\;a_0>\sum\limits_{i=1}^{n}a_i,\;e_0>0,
		\end{align}
		the following inequalities can not hold simultaneously.
		\begin{align}
			\label{lemfinal1}&\dfrac{a_0}{e_0}
			\leq\sum\limits_{i=1}^{n}\dfrac{a_{i}}{e_{i}}-1,\\
			\label{lemfinal2}&\dfrac{a_0}{e_0^2}
			\geq\sum\limits_{i=1}^{n}\dfrac{a_{i}}{e_{i}^2}-1,\;\text{and}\\
			\label{lemfinal3}&\dfrac{a_0}{e_0^3}
			\leq\sum\limits_{i=1}^{n}\dfrac{a_{i}}{e_{i}^3}-1.
		\end{align}
	\end{lemma}
	\begin{lemma}\label{e2<}
		Define \begin{align}
			E(x,\;y,\;z) \;:=\; (1-x)^2(yz-x)(y-z)^2+(1-y)^2(xz-y)(x-z)^2+(1-z)^2(xy-z)(x-y)^2.
		\end{align}Then, for any $x,y,z \in (0,\;1)$,
  we have $E(x,\;y,\;z)<0$.
	\end{lemma}

	Now, we are prepared to prove Theorem \ref{thm:main} (b) (3). According to the hypothesis of Theorem \ref{thm:main} (b) (3), we assume that $S_1$, $S_2$ and $S_3$ are non-empty.
	By \eqref{eq:4}, we have
	\begin{align}
		\label{eq:gg_1231}g(z)\;&= \;  \sum\limits_{i \in S_1}a_i\ln(z+d_i)
		-\sum\limits_{i \in S_2}a_i\ln(-z+d_i)
		+\sum\limits_{i \in S_3}a_i\ln(-z+d_i)
		,\\
		\label{eq:gg_1232}\dfrac{dg}{dz}(z)\;&=\;\sum\limits_{i \in S_1}\dfrac{a_i}{z+d_i}
		+\sum\limits_{i \in S_2}\dfrac{a_i}{-z+d_i}
		-\sum\limits_{i \in S_3}\dfrac{a_i}{-z+d_i}
		.\end{align}Then we have
	\begin{align}\label{d2g123}
		\dfrac{d^2g}{dz^2}(z)\;=\;-\sum\limits_{i \in S_1}\dfrac{a_i}{(z+d_i)^2}
		+\sum\limits_{i \in S_2}\dfrac{a_i}{(-z+d_i)^2}
		-\sum\limits_{i \in S_3}\dfrac{a_i}{(-z+d_i)^2}.
	\end{align}
	
	$\Rightarrow)$First, we prove the sufficiency.  Assume that  there exists a subset $S_2^{*}$ of $S_2$, such that
	\begin{align}\label{con123}
		\sum\limits_{i\in S_{3}}a_i>\sum\limits_{i\in S_2^{*}}a_i>\min\limits_{i\in S_{3}}\{a_i\}.
	\end{align}
 The goal is to prove that 
 $G$ admits multistability. 
	Assume that $a_{p} = \min\limits_{i \in S_3}\{a_i\}$, where
 $p\in S_3$. 
 First, we let $d_p=1$, and for any $i\in S_1$, we let $d_i = 0$. Then,  for any $i \in S_2^*$, we let $d_i = w_1$ and for any $i \in S_2\setminus{S_2^*}$, we let $d_i = w_2$ for any $i \in S_2\setminus{S_2^*}$. Similarly, for any $i \in S_3\setminus{\{p\}}$, we also make all $d_i$'s the same, i.e., we let $d_i = w_3$ for any $i \in S_3\setminus{\{p\}}$. Here, we assume that $w_i>1$ $(i = 1,\;2,\;3)$.
	Then, by \eqref{eq:gg_1232}, we have $$ \dfrac{dg}{dz}(z)\;=\;\dfrac{\sum\limits_{i \in S_1}a_i}{z}
	+\dfrac{\sum\limits_{i \in S_2^{*}}a_i}{-z+w_1}
	+\dfrac{\sum\limits_{i \in S_2\setminus{S_2^{*}}}a_i}{-z+w_2}
	-\dfrac{a_{p}}{-z+1}
	-\dfrac{\sum\limits_{i \in S_3\setminus{\{p\}}}a_i}{-z+w_3}
	.$$
	Note that the interval $I$ defined in \eqref{eq:5} is $(0,\;1)$. By \eqref{eq:gg_1231} and \eqref{eq:gg_1232}, we have $\lim\limits_{z \to 0^+}g(z) = -\infty$, $\lim\limits_{z \to 1^-}g(z) = -\infty$, $\lim\limits_{z \to 0^+}\dfrac{dg}{dz}(z) = +\infty$, and $\lim\limits_{z \to 1^-}\dfrac{dg}{dz}(z) = -\infty$. 
	By Lemma \ref{diff di} and Lemma 4.2 (\ref{lemma4.2_2}), we only need to choose proper numbers $w_1$, $w_2$, and $w_3$ such that there exist $\widetilde{z}_1,\; \widetilde{z}_2\in I$ $(\widetilde{z}_1<\widetilde{z}_2)$ satisfying $\dfrac{dg}{dz}(\widetilde{z}_1)<0$ and $\dfrac{dg}{dz}(\widetilde{z}_2)>0$. 
	Let $$h(z)\;:=\;\dfrac{\sum\limits_{i \in S_1}a_i}{z}
	+\dfrac{\sum\limits_{i \in S_2^{*}}a_i}{-z+w_1}
	-\dfrac{a_{p}}{-z+1}
	-\dfrac{\sum\limits_{i \in S_3\setminus{\{p\}}}a_i}{-z+w_3}
	.$$ We complete the proof by the following three steps. 
	
	(Step 1) In this step, we prove that there exist $w_3>1$ and $\widetilde{z}_1 \in I$ such that for any $w_1>1$, we have $h(\widetilde{z}_1) < 0$. 
	In fact, for any $w_1>1$,  we have 
	\begin{align}\label{hz123}
		h(z)\;<\;\dfrac{\sum\limits_{i \in S_1}a_i}{z}
		+\dfrac{\sum\limits_{i \in S_2^{*}}a_i -a_{p}}{-z+1}
		-\dfrac{\sum\limits_{i \in S_3\setminus{\{p\}}}a_i}{-z+w_3}
		.
	\end{align}
	Let the RHS of \eqref{hz123} be $H(z)$. We solve $d$ from $H(z)<0$, and we get
	\begin{align}\label{d3_123}
		w_3 < \dfrac{\mathcal{N}(z)}{\mathcal{D}(z)},
	\end{align}
 where $\mathcal{N}(z):=\sum\limits_{i \in S_3\setminus{\{p\}}}a_i+\sum\limits_{i \in S_1}a_i+(\sum\limits_{i \in S_2^{*}}a_i -a_{p})\dfrac{z}{-z+1}$ and $\mathcal{D}(z):=\sum\limits_{i \in S_1}a_i\dfrac{1}{z}+(\sum\limits_{i \in S_2^{*}}a_i -a_{p})\dfrac{1}{-z+1}$. 
	Note that 
	\begin{align*}
		&~~~\mathcal{N}(z)- \mathcal{D}(z) \\ 
		&= \sum\limits_{i \in S_3}a_i-\sum\limits_{i \in S_2^{*}}a_i+\sum\limits_{i \in S_1}a_i(1-\dfrac{1}{z}).\;
	\end{align*}
	By \eqref{con123}, we have $\lim\limits_{z\to 1}(\mathcal{N}(z)- \mathcal{D}(z)) = \sum\limits_{i \in S_3}a_i-\sum\limits_{i \in S_2^{*}}a_i>0$. So, there exists $\widetilde{z}_1 \in (0,\;1)$ such that $\mathcal{N}(\widetilde{z}_1)- \mathcal{D}(\widetilde{z}_1) > 0.$
Note also by \eqref{con123}, we have for any $z\in(0,1)$, $\mathcal{N}(z)>0$ and $\mathcal{D}(z)>0$. 
	Hence, $$\dfrac{\mathcal{N}(\widetilde{z}_1)}{\mathcal{D}(\widetilde{z}_1)} > 1.$$
	By \eqref{d3_123}, there exists $w_3>1$ such that $H(\widetilde{z}_1)<0$. Recall that for any $w_1>1$, we have \eqref{hz123}, i.e., $h(z)<H(z)$. So, for any $w_1>1$, we have $h(\widetilde{z}_1)<0$.
	
	(Step 2) In this step, we prove that there exist $w_1>1$ and $\widetilde{z}_2 \in (\widetilde{z}_1,\;1)$ such that $h(\widetilde{z}_2) > 0$.
 In fact, we can solve $w_1$ from $h(z)>0$, and we get
	\begin{align}\label{d1_123_1}
		w_1 < \dfrac{\widetilde{\mathcal{N}}(z)}{\widetilde{\mathcal{D}}(z)}, 
	\end{align}
 where $$\widetilde{\mathcal{D}}(z):=-\sum\limits_{i \in S_1}a_i\dfrac{1}{z}
			+\dfrac{a_{p}}{-z+1}
			+(\sum\limits_{i \in S_3\setminus{\{p\}}}a_i)\dfrac{1}{-z+w_3}$$ and $$\widetilde{\mathcal{N}}(z):=\sum\limits_{i \in S_2^{*}}a_i
			-\sum\limits_{i \in S_1}a_i
			+a_{p}\dfrac{z}{-z+1}
			+(\sum\limits_{i \in S_3\setminus{\{p\}}}a_i)\dfrac{z}{-z+w_3}.$$
	Since 
	\begin{align}
		\lim\limits_{z\to 1^-}\widetilde{\mathcal{D}}(z)=+\infty \;\text{and}\;
  \lim\limits_{z\to 1^-}\widetilde{\mathcal{N}}(z)=+\infty,
	\end{align} there exists $z^*\in (0,\;1)$ such that for any $z \in (z^*,\;1)$,
	\begin{align}
		\widetilde{\mathcal{D}}(z)>0,\;\text{and}\;\widetilde{\mathcal{N}}(z)>0.
	\end{align} 
	 Note that
	\begin{align*}
		~~~&\widetilde{\mathcal{N}}(z)-\widetilde{\mathcal{D}}(z)
		\\=&\sum\limits_{i \in S_2^{*}}a_i-a_{p}+\sum\limits_{i \in S_1}a_i(\dfrac{1}{z}-1)+(\sum\limits_{i \in S_3\setminus{\{p\}}}a_i)\dfrac{z-1}{-z+w_3}.
	\end{align*}
	Since $w_3>1$, we have $$\lim\limits_{z \to 1} (\widetilde{\mathcal{N}}(z)-\widetilde{\mathcal{D}}(z)) = \sum\limits_{i \in S_2^{*}}a_i-a_{p}.$$
	By \eqref{con123}, the above limit is positive. 
	Therefore, we can choose $\widetilde{z}_2 \in (\max\{\widetilde{z}_1,z^*\},\;1)$ such that $$\dfrac{\widetilde{\mathcal{N}}(\widetilde{z}_2)}{\widetilde{\mathcal{D}}(\widetilde{z}_2)} > 1.$$
	 So, by  \eqref{d1_123_1}, we can choose appropriate $w_1 > 1$ such that $h(\widetilde{z}_2) > 0$.
	
	(Step 3) In this step, we prove that there exists $w_2>1$ such that $\dfrac{dg}{dz}(\widetilde{z}_1) < 0$ and $\dfrac{dg}{dz}(\widetilde{z}_2) > 0$. In fact, we can choose $$w_2 = \max{\{\dfrac{2\sum\limits_{i \in S_2\setminus{S_2^{*}}}a_i}{-h(\widetilde{z}_1)}+\widetilde{z}_1,\;2\}}.$$
	Therefore, we have $$\dfrac{dg}{dz}(\widetilde{z}_1) = h(\widetilde{z}_1)+\dfrac{\sum\limits_{i \in S_2\setminus{S_2^{*}}}a_i}{-\widetilde{z}_1+w_2} \leq \dfrac{h(\widetilde{z}_1)}{2} < 0,\;\text{and}$$
	$$\dfrac{dg}{dz}(\widetilde{z}_2)=h(\widetilde{z}_2)+\dfrac{\sum\limits_{i \in S_2\setminus{S_2^{*}}}a_i}{-\widetilde{z}_2+w_2} > h(\widetilde{z}_2) > 0.$$
	
	$\Leftarrow)$ Next, we prove the necessity. Our goal is to prove that if $G$ admits multistability, then there exists a subset $S_2^{*}$ of $S_2$ such that
	\begin{align}\label{s2'}
		\sum\limits_{i\in S_{3}}a_i>\sum\limits_{i\in S_2^{*}}a_i>\min\limits_{i\in S_{3}}\{a_i\}.
	\end{align}
 Assume that $|S_i| = s_i$ 
 $(i=1, 2, 3)$, and assume that $S_1 = \{1,\;\cdots,\;s_1\}$, $S_2 = \{s_1+1,\;\cdots,\;s_1+s_2\}$, and $S_3=\{s_1+s_2+1,\;\cdots,\;s_1+s_2+s_3\}$.
 Below, 
 we prove the conclusion by deducing a contradiction.
	Note that if there does not exist a subset $S_2^{*}$ of $S_2$ such that  $\sum\limits_{i\in S_{3}}a_i>\sum\limits_{i\in S_2^{*}}a_i>\min\limits_{i\in S_{3}}\{a_i\}$, then we have the following three cases.
	\begin{itemize}\label{nosubset}
		\item[(Case 1)] $s_3 = 1.$
    \item[(Case 2)] $s_3\geq 2$ and for any $i\in S_2$, we have $\sum\limits_{i \in S_3}a_i \leq a_i$.
		\item[(Case 3)]  Assume that $a_{s_1+1} \leq a_{s_1+2} \leq \cdots \leq a_{s_1+s_2}$. There exists $ k \in \{1,\;...,\;s_2\}$, such that $\sum\limits_{i=s_1+1}^{s_1+k}a_i \leq \min\limits_{i \in S_3}\{a_i\} < \sum\limits_{i \in S_3}a_i \leq a_{s_1+k+1} \leq \cdots \leq a_{s_1+s_2}.$
	\end{itemize}
  Below, we will prove the conclusion by discussing the three cases. 
 By Lemma \ref{diff di}, if $G$ admits multistability, then there exist $\{d_i\}^s_{i=1}\subset \mathbb{R}$ and $K\in\mathbb{R}$ such that the equation $g(z)=K$ has at least $2$ solutions $z_1\ \text{and}\ z_2$ in the interval $I=(\mathcal{L},\; \mathcal{R})$ defined in \eqref{eq:5} satisfying $\dfrac{dg}{dz}(z_1)<0\ \text{and}\ \dfrac{dg}{dz}(z_2)<0$, where these $d_i$'s are distinct from each other. 
	\begin{enumerate}
		\item[(Case 1)]
		Assume that $s_3 = 1.$ Then, $S_3=\{s_1+s_2+1\}$.
		Suppose $d_{1} < d_{2} < \cdots < d_{s_1}$, and $d_{s_1+1} < d_{s_1+2} < \cdots < d_{s_1+s_2}$.
		\begin{enumerate}
			\item[(Case 1.1)]
			If $d_{s_1+s_2+1} < d_{s_1+1}$, then the interval  $I$ defined in \eqref{eq:5} is  $(-d_{1},\;d_{s_1+s_2+1})$. Notice that by \eqref{eq:gg_1232}, for any $i \in \{s_1+1,\;...,\;s_1+s_2-1\}$, $\lim\limits_{z \to d_i^{+}}{\dfrac{dg}{dz}(z)} = -\infty$ and $\lim\limits_{z \to d_{i+1}^{-}}{\dfrac{dg}{dz}(z)} = +\infty$. Note also that $\dfrac{dg}{dz}(z)$ is continuous in $(d_i, \; d_{i+1})$. So, there exists $ z_i \in (d_i, \; d_{i+1})$ such that $\dfrac{dg}{dz}(z_i) = 0$. Hence, $\dfrac{dg}{dz}(z) = 0$ has at least $s_2-1$ solutions in $(d_{s_1+1},\;+\infty)$. Similarly, notice that by \eqref{eq:gg_1232}, for any $ i \in \{1,\;...,\;s_1-1\}$, $\lim\limits_{z \to -d_i^{-}}{\dfrac{dg}{dz}(z)} = -\infty$ and $\lim\limits_{z \to -d_{i+1}^{+}}{\dfrac{dg}{dz}(z)} = +\infty$.  So, there exists $ z_i \in (-d_{i+1},\;-d_i)$ such that $\dfrac{dg}{dz}(z_i) = 0$. Hence, $\dfrac{dg}{dz}(z) = 0$ has at least $s_1-1$ solutions in $(-\infty,\;-d_{1})$. 
			Since the numerator of $\dfrac{dg}{dz}(z)$ is a polynomial with degree $s_1+s_2$, $\dfrac{dg}{dz}(z)=0$ has no more than $s_1+s_2$ real solutions in $(-\infty,\;+\infty)$.
			 Hence, there are no more than 2 solutions in $I=(-d_{1},\;d_{s_1+s_2+1})$. On the other hand, by \eqref{eq:gg_1231} and \eqref{eq:gg_1232}, notice that $\lim\limits_{z \to -d_{1}^{+}}{g(z)} = -\infty$,  $\lim\limits_{z \to d_{s_1+s_2+1}^{-}}{g(z)} = -\infty$, $\lim\limits_{z \to -d_{1}^{+}}{\dfrac{dg}{dz}(z)} = +\infty$, and $\lim\limits_{z \to d_{s_1+s_2+1}^{-}}{\dfrac{dg}{dz}(z)} = -\infty$. By Lemma \ref{diff di} and Lemma 4.2 (\ref{lemma4.2_2}), if $G$ admits multistability, then there exist $\widetilde{z}_1,\; \widetilde{z}_2\in I$ $(\widetilde{z}_1<\widetilde{z}_2)$ such that $$\dfrac{dg}{dz}(\widetilde{z}_1)<0,\; \text{and}\;\dfrac{dg}{dz}(\widetilde{z}_2)>0.$$ Since $\lim\limits_{z \to -d_{1}^+}\dfrac{dg}{dz}(z) = +\infty$ and $\lim\limits_{z \to d_{s_1+s_2+1}^-}\dfrac{dg}{dz}(z) = -\infty$, $\dfrac{dg}{dz}(z)=0$ has at least $3$ solutions in $I$, which is a contradiction.
			\item[(Case 1.2)]
			If $d_{s_1+s_2+1} > d_{s_1+1}$, then the interval  $I$ defined in \eqref{eq:5} is  $(-d_{1},\;d_{s_1+1})$. Notice that by \eqref{eq:gg_1232}, for any $i \in \{s_1+1,\;...,\;s_1+s_2-1\}$, we have $\lim\limits_{z \to d_i^{+}}{\dfrac{dg}{dz}(z)} = -\infty$ and $\lim\limits_{z \to d_{ i+1}^{-}}{\dfrac{dg}{dz}(z)} = +\infty$. So, for any $i \in \{s_1+1,\;...,\;s_1+s_2-1\}$ satisfying $d_{s_1+s_2+1} \notin (d_i,\;d_{i+1})$, there exists $ z_i \in (d_i,\;d_{i+1})$ such that $\dfrac{dg}{dz}(z_i) = 0$.  Note that $d_{s_1+s_2+1}$ is located in at most  one of the $s_2-1$ intervals $(d_i,\;d_{i+1})\;(i \in \{s_1+1,\;...,\;s_1+s_2-1\})$. Hence, $\dfrac{dg}{dz}(z) = 0$ has at least $s_2-2$ real solutions in $(d_{s_1+1},\;+\infty)$. Similarly, notice that by \eqref{eq:gg_1232}, for any $ i \in \{1,\;...,\;s_1-1\}$, $\lim\limits_{z \to -d_i^{-}}{\dfrac{dg}{dz}(z)} = -\infty$ and $\lim\limits_{z \to -d_{i+1}^{+}}{\dfrac{dg}{dz}(z)} = +\infty$.  So, there exists $ z_i \in (-d_{i+1},\;-d_i)$ such that $\dfrac{dg}{dz}(z_i) = 0$. Hence, $\dfrac{dg}{dz}(z) = 0$ has at least $s_1-1$ real solutions in $(-\infty,\;-d_{1})$. 
			Since the numerator of $\dfrac{dg}{dz}(z)$ is a polynomial with degree $s_1+s_2$, $\dfrac{dg}{dz}(z)=0$ has no more than $s_1+s_2$ real solutions in $(-\infty,\;+\infty)$. Hence, $\dfrac{dg}{dz}(z)=0$ has no more than $3$ real solutions in  $I$. On the other hand, by \eqref{eq:gg_1231} and \eqref{eq:gg_1232}, notice that $\lim\limits_{z \to -d_1^{+}}g(z)=-\infty,\; \lim\limits_{z \to d_{s_1+1}^{-}}g(z)=+\infty,\; \lim\limits_{z \to -d_1^{+}}\dfrac{dg}{dz}(z)=+\infty,\;\text{and} \lim\limits_{z \to d_{s_1+1}^{-}}\dfrac{dg}{dz}(z)=+\infty$. By Lemma \ref{diff di} and Lemma 4.2 (\ref{lemma4.2_3}), if $G$ admits multistability, then there exist $\widetilde{z}_1,\; \widetilde{z}_2,\; \widetilde{z}_3\in I$ $(\widetilde{z}_1<\widetilde{z}_2<\widetilde{z}_3)$ such that $$\dfrac{dg}{dz}(\widetilde{z}_1)<0,\; \dfrac{dg}{dz}(\widetilde{z}_2)>0,\; \text{and}\;\dfrac{dg}{dz}(\widetilde{z}_3)<0.$$ 
   Since $\lim\limits_{z \to -d_{1}^+}\dfrac{dg}{dz}(z) = +\infty$ and $\lim\limits_{z \to d_{s_1+1}^-}\dfrac{dg}{dz}(z) = +\infty$, $\dfrac{dg}{dz}(z)=0$ has at least $4$ solutions in $I$, which is a contradiction.
		\end{enumerate}
\item[(Case 2)]
Recall that the hypothesis of this case is 
   that $s_3\geq 2$ and for any $i\in S_2$, we have $\sum\limits_{i \in S_3}a_i \leq a_i$.
Notice that the interval $I$ defined in \eqref{eq:5} is \begin{align}\label{interval}
			I = (\mathcal{L},\;\mathcal{R}),
		\end{align} where 
		\begin{align}
			&\label{defdm123}\mathcal{L} = -\min\limits_{i \in S_1}\{d_i\},\;\text{and}\\ &\label{defdn123}\mathcal{R} = \min\{d_i\}_{i\in S_2\cup S_3}.
		\end{align}
  (Step 1)
  Below we prove that if $G$ admits multistability, then there exist $\widetilde{z}_1,\; \widetilde{z}_2,\; \widetilde{z}_3\in I$ $(\widetilde{z}_1<\widetilde{z}_2<\widetilde{z}_3)$, such that
		\begin{align}\label{<><}
			\dfrac{dg}{dz}(\widetilde{z}_1)<0,\; \dfrac{dg}{dz}(\widetilde{z}_2)>0,\;\text{and}\; \dfrac{dg}{dz}(\widetilde{z}_3)<0.
		\end{align} 
 If $\min\limits_{i \in S_2}\{d_i\}<\min\limits_{i \in S_3}\{d_i\}$, we have $I = (-\min\limits_{i \in S_1}\{d_i\},\;\min\limits_{i \in S_2}\{d_i\})$. By \eqref{eq:gg_1231} and \eqref{eq:gg_1232}, we have $\lim\limits_{z \to \mathcal{L}^+}g(z) = -\infty$, $\lim\limits_{z \to \mathcal{R}^-}g(z) = +\infty$, $\lim\limits_{z \to \mathcal{L}^+}\dfrac{dg}{dz}(z) = +\infty$, and $\lim\limits_{z \to \mathcal{R}^-}\dfrac{dg}{dz}(z) = +\infty$. By Lemma 4.2 (\ref{lemma4.2_3}), 
there exist $\widetilde{z}_1,\; \widetilde{z}_2,\; \widetilde{z}_3\in I$ $(\widetilde{z}_1<\widetilde{z}_2<\widetilde{z}_3)$, such that
		\begin{align}
			\dfrac{dg}{dz}(\widetilde{z}_1)<0,\; \dfrac{dg}{dz}(\widetilde{z}_2)>0,\;\text{and}\; \dfrac{dg}{dz}(\widetilde{z}_3)<0.
		\end{align} 
 If $\min\limits_{i \in S_3}\{d_i\}<\min\limits_{i \in S_2}\{d_i\}$, we have $I = (-\min\limits_{i \in S_1}\{d_i\},\;\min\limits_{i \in S_3}\{d_i\})$. By \eqref{eq:gg_1231} and \eqref{eq:gg_1232}, we have $\lim\limits_{z \to \mathcal{L}^+}g(z) = -\infty$, $\lim\limits_{z \to \mathcal{R}^-}g(z) = -\infty$, $\lim\limits_{z \to \mathcal{L}^+}\dfrac{dg}{dz}(z) = +\infty$, and $\lim\limits_{z \to \mathcal{R}^-}\dfrac{dg}{dz}(z) = -\infty$. By Lemma 4.2 (\ref{lemma4.2_2}), 
there exist $\widetilde{z}_1,\; \widetilde{z}_2\in I$ $(\widetilde{z}_1<\widetilde{z}_2)$, such that
		\begin{align}
			\dfrac{dg}{dz}(\widetilde{z}_1)<0,\; \text{and}\;\dfrac{dg}{dz}(\widetilde{z}_2)>0.
		\end{align} 
  Since $\lim\limits_{z \to \mathcal{R}^-}\dfrac{dg}{dz}(z) = -\infty$, there exists $\widetilde{z}_3\in (\widetilde{z}_2,\;R)$, such that $\dfrac{dg}{dz}(\widetilde{z}_3)<0$.
So, the conclusion holds.
  
  (Step 2) Below, we deduce a contradiction by \eqref{<><}. 
Note that $\dfrac{dg}{dz}(z)$ is a rational function. So, by \eqref{<><}, there exists $z_0\in(\widetilde{z}_1,\;\widetilde{z}_2)$, such that
			\begin{align}\label{dg=>}
				\dfrac{dg}{dz}(z_0)=0,\;\text{and}\dfrac{d^2g}{dz^2}(z_0)\geq0.
			\end{align} 
   By \eqref{eq:gg_1232}, \eqref{d2g123} and \eqref{dg=>}, we have 
   \begin{align}
       \sum\limits_{i \in S_1}\dfrac{a_i}{z_0+d_i}
		+\sum\limits_{i \in S_2}\dfrac{a_i}{-z_0+d_i}
		-\sum\limits_{i \in S_3}\dfrac{a_i}{-z_0+d_i}=0,\;\text{and}\\
  -\sum\limits_{i \in S_1}\dfrac{a_i}{(z_0+d_i)^2}
		+\sum\limits_{i \in S_2}\dfrac{a_i}{(-z_0+d_i)^2}
		-\sum\limits_{i \in S_3}\dfrac{a_i}{(-z_0+d_i)^2}\geq0.
   \end{align}
   Recall that $a_i>0$ and $z_0+d_i>0\;(i\in S_1)$. Then, we have
      \begin{align}
      \label{more1} \sum\limits_{i \in S_2}\dfrac{a_i}{-z_0+d_i}
		<\sum\limits_{i \in S_3}\dfrac{a_i}{-z_0+d_i},\;\text{and}\\
\label{more2}\sum\limits_{i \in S_2}\dfrac{a_i}{(-z_0+d_i)^2}
		>\sum\limits_{i \in S_3}\dfrac{a_i}{(-z_0+d_i)^2}.
   \end{align}
Recall that the hypothesis of this case is that $s_3\geq2$ and for any $j \in S_2$, we have $a_j\geq\sum\limits_{i\in S_3}a_i$. Then, by \eqref{more1}--\eqref{more2} and by  Cauchy's inequality, we have \begin{align}
    \sum\limits_{i \in S_2}\dfrac{a_i^2}{(-z_0+d_i)^2}\geq \sum\limits_{i \in S_2}\dfrac{a_i}{(-z_0+d_i)^2}\sum\limits_{i\in S_3}a_i>\sum\limits_{i \in S_3}\dfrac{a_i}{(-z_0+d_i)^2}\sum\limits_{i\in S_3}a_i\\ \geq
    (\sum\limits_{i \in S_3}\dfrac{a_i}{-z_0+d_i})^2>(\sum\limits_{i \in S_2}\dfrac{a_i}{-z_0+d_i})^2,
\end{align}
which is a contradiction.
		\item[(Case 3)]
  Recall that the hypothesis of this case is 
   that  there exists $ k \in \{1,\;...,\;s_2\}$, such that
		\begin{align}\label{conleq}
			\sum\limits_{i=s_1+1}^{s_1+k}a_i \leq \min\limits_{i \in S_3}\{a_i\} < \sum\limits_{i \in S_3}a_i \leq a_{s_1+k+1} \leq \cdots \leq a_{s_1+s_2},
		\end{align} 
  where  $a_{s_1+1} \leq a_{s_1+2} \leq \cdots \leq a_{s_1+s_2}$.
  Similar to the proof of Case 2, we have \eqref{<><}.
  
		(Step 1) In this step, for the $\widetilde{z}_1,\;\widetilde{z}_2,\;\widetilde{z}_3$ in \eqref{<><}, we prove that there exists \begin{align}\label{G_6'}
				\widetilde{G}(z) \;:=\;\dfrac{A_1}{z+D_1}
				+\dfrac{A_2}{-z+D_2}
				+\dfrac{A_3}{-z+D_3}+\theta
				-\sum\limits_{i \in S_3}\dfrac{a_i}{-z+d_i}, 
			\end{align}
   such that \begin{align}
				\label{gz><}\widetilde{G}(\widetilde{z}_1)<0,\;\widetilde{G}(\widetilde{z}_2)>0,\;\widetilde{G}(\widetilde{z}_3)<0,
    	\end{align}
   where
			\begin{align}
&\label{G1}A_1>0,\;0<A_2\leq\min\limits_{i \in S_3}\{a_i\},\;A_3\geq\sum\limits_{i \in S_3}a_i,\\
    &\label{G2}\widetilde{z}_i+D_1>0,\;-\widetilde{z}_i+D_2>0,\;-\widetilde{z}_i+D_3>0\;\;(i=1,\;2,\;3),\;\text{and}\\
				&\label{G3}\theta\geq0.
			\end{align}

		
		Recall that by \eqref{eq:ar}, for any $  i \in S_1$, we have $a_{i}>0$. Since $\widetilde{z}_i\in I$ ($i=1,\;2,\;3$), by \eqref{defdm123}, we have $\widetilde{z}_i>-\min\limits_{i\in S_1}\{d_{i}\}$ ($i=1,\;2,\;3$). Then, by Lemma \ref{simplify3}, there exist $A_1,\;D_1,\;\theta_1\in R$ such that \begin{align}\label{change1}
			\sum\limits_{i \in S_1}\dfrac{a_i}{\widetilde{z}_j+d_i} = \dfrac{A_1}{\widetilde{z}_j+D_1}+\theta_1 \;\; (j = 1,\;2,\;3),
		\end{align} 
		where
		\begin{align}
			\label{change1'}
   \min\limits_{i\in{S_{1}}}\{a_i\}\leq A_1\leq\sum\limits_{i\in{S_{1}}}a_i
   ,\;D_1\geq\min\limits_{i\in S_1}\{d_{i}\},\;\text{and}\;\theta_1\geq0.
		\end{align}
 So, we have $A_1>0$ in \eqref{G1} and $\widetilde{z}_i+D_1>0\;\;(i = 1,\;2,\;3)$ in \eqref{G2}.  
		By \eqref{conleq}, let $S_{2}^{(1)} = \{s_1+1,\;\cdots,\;s_1+k\}$ and $S_{2}^{(2)} = \{s_1+k+1,\;\cdots,\;s_1+s_2\}$. So, $S_2 = S_{2}^{(1)} \cup S_{2}^{(2)}$ and by \eqref{conleq}, we have
		\begin{align}
			\label{S2a}\sum\limits_{i\in S_{2}^{(1)}}a_i \leq \min\limits_{i \in S_3}\{a_i\}, \\ 
			\label{S2b}\sum\limits_{i \in S_3}a_i \leq \min\limits_{i \in S_{2}^{(2)}}\{a_i\}.
		\end{align}\\
		Recall that by \eqref{eq:ar}, for any $  i \in S_2$, $a_i>0$. Since $\widetilde{z}_i\in I$ ($i=1,\;2,\;3$), by \eqref{defdn123},  we have $-\widetilde{z}_j>-\mathcal{R}\geq-\min\limits_{i\in S_{2}^{(1)}}\{d_{i}\}$ ($j=1,\;2,\;3$). Similarly, by Lemma \ref{simplify3}, there exist $A_2,\;D_2,\;\theta_2\in R$ such that 
		\begin{align}\label{change2}
			\sum\limits_{i \in S_{2}^{(1)}}\dfrac{a_i}{-\widetilde{z}_j+d_i} = \dfrac{A_2}{-\widetilde{z}_j+D_2}+\theta_2 \;\;( j = 1,\;2,\;3),
		\end{align}
		where
		\begin{align}
			\label{change2'}\min\limits_{i\in{S_{2}^{(1)}}}\{a_i\}\leq A_2\leq\sum\limits_{i\in{S_{2}^{(1)}}}a_i,\;D_2\geq\min\limits_{i\in S_{2}^{(1)}}\{d_{i}\},\;\text{and}\;\theta_2\geq0.
		\end{align}
  Then, by \eqref{S2a}, we have $0<A_2\leq\min\limits_{i \in S_3}\{a_i\}$ in \eqref{G1} and $-\widetilde{z}_i+D_2>0\;\;(i = 1,\;2,\;3)$ in \eqref{G2}.
		Recall that by \eqref{eq:ar}, for any $  i \in S_2$, $a_i>0$. Since $\widetilde{z}_i\in I$ ($i=1,\;2,\;3$), by \eqref{defdn123},  we have $-\widetilde{z}_j>-\mathcal{R}\geq-\min\limits_{i\in S_{2}^{(2)}}\{d_{i}\}$ ($j=1,\;2,\;3$). Similarly, by Lemma \ref{simplify3}, there exist $A_3,\;D_3,\;\theta_3\in R$ such that 
		\begin{align}\label{change3}
			\sum\limits_{i \in S_{2}^{(2)}}\dfrac{a_i}{-\widetilde{z}_j+d_i} = \dfrac{A_3}{-\widetilde{z}_j+D_3}+\theta_3 \;\;(j = 1,\;2,\;3),
		\end{align}
		where
		\begin{align}
			\label{change3'}\min\limits_{i\in{S_{2}^{(2)}}}\{a_i\}\leq A_3\leq\sum\limits_{i\in{S_{2}^{(2)}}}a_i,\;D_3\geq\min\limits_{i\in S_{2}^{(2)}}\{d_i\},\;\text{and}\;\theta_3\geq0.
		\end{align}
   Then, by \eqref{S2b}, we have $A_3\geq\sum\limits_{i \in S_3}\{a_i\}$ in \eqref{G1} and $-\widetilde{z}_i+D_3>0\;\;(i = 1,\;2,\;3)$ in \eqref{G2}.
		Let 
		$	\theta := \theta_1+\theta_2+\theta_3.$
		By \eqref{change1'}, \eqref{change2'}, and \eqref{change3'},  we have \eqref{G3}.
  By \eqref{change1}, \eqref{change2}, and \eqref{change3}, we have $\widetilde{G}(\widetilde{z}_i) = \dfrac{dg}{dz}(\widetilde{z}_i)$ ($i = 1,\;2,\;3$).
		Then, by \eqref{<><}, we have \eqref{gz><}.
  
 (Step 2) Below we will prove \eqref{gz><} can not hold by discussing two cases. 
  
		\begin{enumerate}
			\item[(Case 2.1)]
			Assume that
			\begin{align}\label{M*=d1}
				\min\limits_{i \in S_{3}}\{d_i\}\leq D_2.
			\end{align} 
			
			(Step 1) In this step, we prove that there exists $z^* \in (\widetilde{z}_1,\;\widetilde{z}_2)$, such that
   \begin{align}
      \label{G_6=0}&\dfrac{A_2}{-z^*+D_2}+
				\dfrac{A_3}{-z^*+D_3}
				-\sum\limits_{i \in S_3}\dfrac{a_i}{-z^*+d_i}<0,\;\text{and}\\
    \label{G_6<0}&\dfrac{A_2}{(-z^*+D_2)^2}+\dfrac{A_3}{(-z^*+D_3)^2}-\sum\limits_{i \in S_3}\dfrac{a_i}{(-z^*+d_i)^2}\geq0.
   \end{align}
			Note that by \eqref{gz><}, we have  $
			\widetilde{G}(\widetilde{z}_1)<0,\;\widetilde{G}(\widetilde{z}_2)>0$. Note that $\widetilde{G}(z)$ is a 
rational function. So, there exists $z^* \in (\widetilde{z}_1,\;\widetilde{z}_2)$, such that 
			\begin{align}
				\label{G=0}&\widetilde{G}(z^*) = 0,\\
				\label{dG>0}&\dfrac{d\widetilde{G}}{dz}(z^*) \geq 0.
			\end{align} 
			Since $z^* \in (\widetilde{z}_1,\;\widetilde{z}_2)$, we have $z^*+D_1>0$. Note that by \eqref{G1}, we have $A_1>0$. By \eqref{G=0}, we have \eqref{G_6=0}.
			Notice that 
			\begin{align}\label{def:dG}
				\dfrac{d\widetilde{G}}{dz} = -\dfrac{A_1}{(z+D_1)^2}+\dfrac{A_2}{(-z+D_2)^2}
				+\dfrac{A_3}{(-z+D_3)^2}
				-\sum\limits_{i \in S_3}\dfrac{a_i}{(-z+d_i)^2}.
			\end{align}
		  So, by \eqref{dG>0} and \eqref{def:dG}, we have \eqref{G_6<0}.
   
   (Step 2) Note that there exists $q \in S_3$, such that $d_{q} = \min\limits_{i \in S_3}\{d_i\}$. Let $S_{3}^{*} = S_{3}\setminus{\{q\}}$. 
			Let $e_{2}^{(1)} = \dfrac{-z^*+d_{q}}{-z^*+D_2}$, $e_{2}^{(2)} = \dfrac{-z^*+d_{q}}{-z^*+D_3}$ and let  $e_{i} = \dfrac{-z^*+d_{q}}{-z^*+d_i} $ for any $ i \in S_3^{*}$.
   In this step we prove that 
\begin{align}\label{eq:threezerosmall}\sum\limits_{i \in S_3}a_i(\sum\limits_{i \in S_3^{*}}a_ie_{i}^2+a_{q}-A_2(e_{2}^{(1)})^{2})
				- (\sum\limits_{i \in S_3^{*}}a_ie_{i}
				+a_{q}-A_2e_{2}^{(1)})^2<0.\end{align}
    
			We multiply \eqref{G_6=0} by $-z+d_{q}$ and we can get
			\begin{align}
				\label{g<123} &
				A_3e_{2}^{(2)}
				<\sum\limits_{i \in S_3^{*}}a_ie_{i}
				+a_{q}-A_2e_{2}^{(1)}.
			\end{align}
   By \eqref{G1}, we have 
			\begin{align}
				\label{a2b123}A_3 \geq \sum\limits_{ i \in S_3}a_i.
			\end{align} By \eqref{g<123} and \eqref{a2b123}, we have $$A_3(e_{2}^{(2)})^{2} = {A_3}^{2}(e_{2}^{(2)})^{2}\dfrac{1}{A_3}
			<(\sum\limits_{i \in S_3^{*}}a_ie_{i}
			+a_{q}-A_2e_{2}^{(1)})^2\dfrac{1}{\sum\limits_{i \in S_3}a_i}.$$
We multiply \eqref{G_6<0} by $(-z+d_{q})^2$ and we can get		
			\begin{align}
				\label{g>123} A_3(e_{2}^{(2)})^2
				\geq\sum\limits_{i \in S_3^{*}}a_ie_{i}^2
				+a_{q}-A_2(e_{2}^{(1)})^2.
			\end{align}
			 Then, we have
			$$\sum\limits_{i \in S_3^{*}}a_ie_{i}^2+a_{q}-A_2(e_{2}^{(1)})^{2}
			<(\sum\limits_{i \in S_3^{*}}a_ie_{i}
			+a_{q}-A_2e_{2}^{(1)})^2\dfrac{1}{\sum\limits_{i \in S_3}a_i}.$$
 Hence, we have \eqref{eq:threezerosmall}.
			
			(Step 2) 
			Define $$\mathcal{F}(y) := (\sum\limits_{i \in S_3^{*}}a_i+y)(\sum\limits_{i \in S_3^{*}}a_ie_{i}^2+y-A_2(e_{2}^{(1)})^{2})
			- (\sum\limits_{i \in S_3^{*}}a_ie_{i}
			+y-A_2e_{2}^{(1)})^2.$$
  Note that $\mathcal{F}(a_{q})<0$ is equivalent to  \eqref{eq:threezerosmall}.  In order to deduce a contradiction,  the goal of this step is to prove $F(a_q)\geq 0$.
    Notice that by \eqref{G1}, we have $a_q\geq A_2$.	 So, we  only need to show that  $\dfrac{d\mathcal{F}}{dy}\geq0$ for any $y\in {\mathbb R}$  and $\mathcal{F}(A_2)\geq0$. 
			Note that 
			\begin{align*}
				\dfrac{d\mathcal{F}}{dy}
				=\sum\limits_{i \in S_3^{*}}a_i(e_{i}-1)^2+A_2e_{2}^{(1)}+A_2e_{2}^{(1)}(1-e_{2}^{(1)}),
			\end{align*}
   and note that by \eqref{M*=d1}, we have 
				$0<e_{2}^{(1)}\leq1$.
			 Note also that by \eqref{G1},  $A_2>0$. So,  $$\dfrac{dF}{dy} > \sum\limits_{i \in S_3^{*}}a_i(e_{i}-1)^2 \geq 0$$
	Below, we prove that ${\mathcal F}(A_2)\geq 0$. 		Define 
			\begin{align}
				&\mathcal{G}(x) := (\sum\limits_{i \in S_3^{*}}a_i+A_2)(\sum\limits_{i \in S_3^{*}}a_ie_{i}^2+A_2-A_2x^{2})
				- (\sum\limits_{i \in S_3^{*}}a_ie_{i}
				+A_2-A_2x)^2. 
			\end{align}
    Notice that $\mathcal{G}(\dfrac{1}{e_{2}^{(1)}}) = \mathcal{F}(A_2)$. So, we only need to prove that $\mathcal{G}(\dfrac{1}{e_{2}^{(1)}})\geq0$.  
			Notice that $\mathcal{G}(x)$ is a quadratic function in $x$, and its coefficient of $x^2$ is negative. So the minimum of $\mathcal{G}$ over $(0,\;1)$ is greater than $\mathcal{G}(0)$ and $\mathcal{G}(1)$.
			By Cauchy's inequality, we have $$\mathcal{G}(0) = (\sum\limits_{i \in S_3^{*}}a_i+A_2)(\sum\limits_{i \in S_3^{*}}a_ie_{i}^2+A_2)
			- (\sum\limits_{i \in S_3^{*}}a_ie_{i}
			+A_2)^2 \geq 0.$$ Also,\; $$\mathcal{G}(1) = (\sum\limits_{i \in S_3^{*}}a_i+A_2)\sum\limits_{i \in S_3^{*}}a_ie_{i}^2
			- (\sum\limits_{i \in S_3^{*}}a_ie_{i})^2>\sum\limits_{i \in S_3^{*}}a_i\sum\limits_{i \in S_3^{*}}a_ie_{i}^2
			- (\sum\limits_{i \in S_3^{*}}a_ie_{i})^2 \geq 0.$$ So, for any $x\in (0,\;1]$, $\mathcal{G}(x)\geq0$. Notice that $\dfrac{1}{e_{2}^{(1)}}\in (0,\;1]$. Then, we have $\mathcal{G}(\dfrac{1}{e_{2}^{(1)}})\geq0$.
			\item[(Case 2.2)]
			Assume that 
			\begin{align}\label{M=d2a}
				\min\limits_{i \in S_{3}}\{d_i\}> D_2.
			\end{align}

   Define 
			\begin{align}\label{pz}
				P(z) := \widetilde{G}(z)-\theta
			\end{align}\\
			Note that  \eqref{gz><} is equivalent to 
			\begin{align}\label{p<><}
				P(\widetilde{z}_1)<-\theta,\;P(\widetilde{z}_2)>-\theta,\;\text{and}\;P(\widetilde{z}_3)<-\theta,
			\end{align} 
   where $\widetilde{z}_1<\widetilde{z}_2<\widetilde{z}_3$.
 Below, we  prove that the three inequalities in  \eqref{p<><} can not hold 
concurrently by discussing two cases.
			\begin{enumerate}
				\item[(Case 2.3.1)] We assume that for any $z\in (\widetilde{z}_1,\;\widetilde{z}_3)$, we have $P(z)\leq0$.
    Note that $P(z)$ is a rational function. So, if \eqref{p<><} holds, then there exists $z_0\in(\widetilde{z}_1,\;\widetilde{z}_3)$, such that
			\begin{align}
				\label{dp=}&\dfrac{dP}{dz}(z_0)=0,\;\text{and}\\
				\label{d2p<}&\dfrac{d^2P}{dz^2}(z_0)\leq0.
			\end{align} 
				Since $z_0\in (\widetilde{z}_1,\;\widetilde{z}_3)$, we have $P(z_0)\leq0$. Then, by \eqref{G_6'} and \eqref{pz}, we have 
				\begin{align}\label{pz<0}
					\dfrac{A_1}{z_0+D_1}+\dfrac{A_2}{-z_0+D_2}
					+\dfrac{A_3}{-z_0+D_3}
					\leq\sum\limits_{i \in S_3}\dfrac{a_i}{-z_0+d_i}.
				\end{align}
    Note that  by \eqref{G1} and \eqref{G2}, we have $\dfrac{A_1}{z_0+D_1}>0$. Then, by \eqref{pz<0}, we have 
				\begin{align}\label{final8}
					\dfrac{A_2}{-z_0+D_2}
					+\dfrac{A_3}{-z_0+D_3}
					\leq\sum\limits_{i \in S_3}\dfrac{a_i}{-z_0+d_i}.
				\end{align}
				Note that $
					\dfrac{dP}{dz} = \dfrac{d\widetilde{G}}{dz}.$
				By \eqref{def:dG} and \eqref{dp=}, we have 
				\begin{align}\label{final9}
					\sum\limits_{i \in S_3}\dfrac{a_i}{(-z_0+d_i)^2}=-\dfrac{A_{1}}{(z_0+D_1)^2}+\dfrac{A_2}{(-z_0+D_2)^2}
					+\dfrac{A_3}{(-z_0+D_3)^2}\leq\dfrac{A_2}{(-z_0+D_2)^2}
					+\dfrac{A_3}{(-z_0+D_3)^2}
					.
				\end{align}
				By \eqref{d2p<}, we have
				\begin{align}
					\dfrac{d^2P}{dz^2}(z_0)=\dfrac{2A_{1}}{(z_0+D_1)^3}+
					\dfrac{2A_2}{(-z_0+D_2)^3}
					+\dfrac{2A_3}{(-z_0+D_3)^3}
					-\sum\limits_{i \in S_3}\dfrac{2a_i}{(-z_0+d_i)^3}\leq 0.
				\end{align}
				Then, we have 
				\begin{align}\label{final10}
					\dfrac{A_2}{(-z_0+D_2)^3}
					+\dfrac{A_3}{(-z_0+D_3)^3}
					&\leq\sum\limits_{i \in S_3}\dfrac{a_i}{(-z_0+d_i)^3}.
				\end{align}
				Below we prove that \eqref{final8}, \eqref{final9}, and \eqref{final10} can not hold concurrently.
				Let $\gamma=\dfrac{A_3}{A_2},\; 
				e = \dfrac{-z_0
					+D_3}{-z_0+D_2}$ and for any $i\in S_3$, let $
     \gamma_{i}=\dfrac{a_i}{A_2},\; 
				e_{i} = \dfrac{-z_0+d_i}{-z_0+D_2}$. By \eqref{G1}, we have $\gamma_i>1\;\;(i\in S_3)$ and $\gamma>\sum\limits_{i\in S_3}\gamma_i$.  By \eqref{G2} and \eqref{M=d2a}, we have $e_i>1\;\;(i\in S_3)$ and $e>0$. Then, we divide the both sides of \eqref{final8} by $\frac{A_2}{-z_0+D_2}$, and we get \begin{align}\label{final6}
					1+\dfrac{\gamma}{e}
					\leq\sum\limits_{i\in S_3}\dfrac{\gamma_{i}}{e_{i}}.
				\end{align}
				Similarly, by \eqref{final9} and \eqref{final10}, we have 
				\begin{align}
					&\label{final7}1+\dfrac{\gamma}{e^2}
					\geq\sum\limits_{i\in S_3}\dfrac{\gamma_{i}}{e_{i}^2},\;\text{and}\\
					&\label{final7'}1+\dfrac{\gamma}{e^3}
					\leq\sum\limits_{i\in S_3}\dfrac{\gamma_{i}}{e_{i}^3}.
				\end{align}
				By Lemma \ref{finalineq}, \eqref{final6}, \eqref{final7} and \eqref{final7'} lead to a contradiction.
				\item[(Case 2.3.2)]We assume that there exists $z_0 \in (\widetilde{z}_1,\;\widetilde{z}_3)$ such that $P(z_0)>0$.
				Let 
				\begin{align}
					\label{hP}\psi(z) &:= (z+D_1)P(z) \\
					\label{defh2}&=A_1+ (z+D_1)(\dfrac{A_2}{-z+D_2}
					+\dfrac{A_3}{-z+D_3}-
					\sum\limits_{i \in S_3}\dfrac{a_i}{-z+d_i}).\end{align}
    Note that $P(z_0)>0$. Recall that by \eqref{G2}, we have $\theta > 0$. So, by \eqref{p<><}, we have 
$ P(\widetilde{z}_1)<0\;\text{and}\;P(\widetilde{z}_3)<0.$ Then by \eqref{hP}, we have 
				\begin{align}\label{H<><}
					&\psi(\widetilde{z}_1)<0,\; \psi(z_0)>0\;\text{and}\; \psi(\widetilde{z}_3)<0.
				\end{align}
				 Note that $\psi(z)$ is a rational function. So, by \eqref{H<><}, there exists $z_0^*\in (\widetilde{z}_1,\;\widetilde{z}_3)$, such that
				\begin{align}
					\label{H>} &	\psi(z_0^*)>0,\\
					\label{dH=} &\dfrac{d\psi}{dz}(z_0^*) = 0,\;\text{and}\\
					\label{d2H<} &\dfrac{d^2\psi}{dz^2}(z_0^*)\leq0.
				\end{align} 
				(Step 1) In this step,  we prove that by \eqref{H<><} and \eqref{H>}, we have 
    \begin{align}\label{h<w}
        \psi(z_0^*)\leq A_1-A_2-A_3+\sum\limits_{i\in S_3}a_i.
    \end{align} We will prove the conclusion by deducing a contradiction.  Let
     \begin{align}\label{defw}
         w &:= A_1-A_2-A_3+\sum\limits_{i\in S_3}a_i.
     \end{align}
				Assume that
				\begin{align}\label{h>w}
					\psi(z_0^*) > w.
				\end{align}
    We first prove that for any $u\in R$,  $\psi(z) = u$ has at most $4$ solutions in $(-\infty,\;D_2)$.
    Notice that $S_3 = \{s_1+s_2+1,\;\cdots,\;s_1+s_2+s_3\}$. Note that by \eqref{diff di}, $d_i$'s are distinct from each other. Assume that $d_{s_1+s_2+1}<d_{s_1+s_2+2}<\cdots<d_{s_1+s_2+s_3}$. Recall that by \eqref{M=d2a}, we have $D_2<d_{s_1+s_2+1}$. Notice that by \eqref{defh2}, for any $i \in S_3$, we have $\lim\limits_{z \to d_i^{+}}{\psi(z)} = -\infty$ and $\lim\limits_{z \to d_{ i+1}^{-}}{\psi(z)} = +\infty$. So, for any $i \in S_3$ satisfying $D_3 \notin (d_i,\;d_{i+1})$, there exists $ z_i \in (d_i,\;d_{i+1})$ such that $\psi(z_i) = u$.  Note that $D_3$ is located in at most  one of the $s_3-1$ intervals $(d_i,\;d_{i+1})$ $(i \in \{s_1+s_2+1,\;\cdots,\;s_1+s_2+s_3-1\})$. Hence, $\psi(z) = u$ has at least $s_3-2$ real solutions in $(D_2,\;+\infty)$.
    Since the numerator of $\psi(z)-u$ is a polynomial with degree $s_3+2$, $\psi(z)=u$ has at most $s_3+2$ solutions in $(-\infty,\;+\infty)$.
    Hence, $\psi(z) = u$ has at most $4$ solutions in $(-\infty,\;D_2)$. Let 
	\begin{align}\label{defu}
					u := \frac{\max\{w,\;\psi(\widetilde{z}_1),\;\psi(\widetilde{z}_3)\}+\min\{A_1,\;\psi(z_0^*)\}}{2}.
				\end{align} 
    Below, we prove that if \eqref{h>w} holds, then $\psi(z) = u$ has at least $5$ solutions in $(-\infty,\;D_2)$, which will be a contradiction.
 We first prove that 
				\begin{align}\label{uinterval}
				u>	\max\{w,\;\psi(\widetilde{z}_1),\;\psi(\widetilde{z}_3)\} ,\;\text{and}\;u<\min\{A_1,\;\psi(z_0^*)\}.
				\end{align} By \eqref{defu}, we only need to prove $\min\{A_1,\;\psi(z_0^*)\}>\max\{w,\;\psi(\widetilde{z}_1),\;\psi(\widetilde{z}_3)\}$.  By \eqref{H<><} and \eqref{H>}, we have $\psi(z_0^*)>\max\{\psi(\widetilde{z}_1),\;\psi(\widetilde{z}_3)\}$. Then, by \eqref{h>w}, we have 
				\begin{align}\label{h>}
					\psi(z_0^*)>\max\{w,\;\psi(\widetilde{z}_1),\;\psi(\widetilde{z}_3)\}.
				\end{align} 
    
    Notice that by \eqref{G1}, we have $A_1>0$. Then by \eqref{H<><}, we have $A_1>\max\{\psi(\widetilde{z}_1),\;\psi(\widetilde{z}_3)\}$. By \eqref{defw}, we have
				\begin{align}\label{A2+A3}
					A_1-w=A_2+A_3-\sum\limits_{i \in S_3}a_i.
				\end{align} By \eqref{G1}, we have $A_2+A_3-\sum\limits_{i \in S_3}a_i>0$. So, we have $A_1>w$.  Then, we have
				\begin{align}\label{A_1>}
					A_1>\max\{w,\;\psi(\widetilde{z}_1),\;\psi(\widetilde{z}_3)\}.
				\end{align} 
				Hence, by \eqref{h>} and \eqref{A_1>}, we have $\min\{A_1,\;\psi(z_0^*)\}>\max\{w,\;\psi(\widetilde{z}_1),\;\psi(\widetilde{z}_3)\}$. Then, \eqref{uinterval} holds.
				Note that by \eqref{uinterval}, we have 
    \begin{align}\label{u1}
u>\psi(\widetilde{z}_1),\;u<\psi(z_0^*),\;u>\psi(\widetilde{z}_3).
    \end{align}
				Note that by \eqref{hP} and \eqref{defw}, we have 
    $\lim\limits_{z\to -\infty}\psi(z) = w$, 
    $\psi(-D_1) = A_1$, and $\lim\limits_{z\to \min\{D_2,\;D_3\}^-}\psi(z) = +\infty$. Then, by \eqref{uinterval}, we have \begin{align*}	u>\lim\limits_{z\to -\infty}\psi(z),\;\text{and}\;u<\psi(-D_1). 
				\end{align*}
    By \eqref{G2}, we have $-\infty<-D_1<\widetilde{z}_1<z_0^*<\widetilde{z}_3<\min\{D_2,\;D_3\}$. Then, we have $\psi(z) = u$ has at least $5$ solutions in $(-\infty,\;\min\{D_2,\;D_3\})\subseteq(-\infty,\;D_2)$, which is a contradiction.
				So, we have \eqref{h<w}.
				
				(Step 2) Let
				\begin{align}\label{defh'}
					\widetilde{\psi}(z)\;:=\;\dfrac{A_2}{-z+D_2}
					+\dfrac{A_3}{-z+D_3}-
					\sum\limits_{i \in S_3}\dfrac{a_i}{-z+d_i}.
				\end{align} In this step,  we prove that by \eqref{h<w}, we have \begin{align}\label{final1}
					\widetilde{\psi}(z_0^*)=\dfrac{A_2}{-z_0^*+D_2}
					+\dfrac{A_3}{-z_0^*+D_3}-
					\sum\limits_{i \in S_3}\dfrac{a_i}{-z_0^*+d_i}<0.
				\end{align}
				By \eqref{defh2}, it is straightforward to check that
				\begin{align}\label{h=}
    \psi(z) = A_1+(z+D_1)\widetilde{\psi}(z).
				\end{align}
				By \eqref{h<w} and \eqref{h=}, we have 
				\begin{align}\label{H<*}
					A_2+A_3-\sum\limits_{i\in S_3}a_i+(z_0^*+D_1)\widetilde{\psi}(z_0^*)\leq0.
				\end{align} 
				Since $z_0^*\in (\widetilde{z}_1,\;\widetilde{z}_3)$, by \eqref{G2}, we have $z_0^*+D_1>0$. By \eqref{G1}, we have $A_2+A_3-\sum\limits_{i\in S_3}a_i>0$. Then, by \eqref{H<*}, we have 
				\eqref{final1}.
    
    (Step 3)
In this step, we prove that by \eqref{dH=}, \eqref{d2H<}, and
 \eqref{final1}, we have
 \begin{align}
     \label{final2}&\frac{d\widetilde{\psi}}{dz}(z_0^*)=\dfrac{A_2}{(-z_0^*+D_2)^2}
					+\dfrac{A_3}{(-z_0^*+D_3)^2}-\sum\limits_{i \in S_3}\dfrac{a_i}{(-z_0^*+d_i)^2}>0,\;\text{and}\\
     \label{final3}&\frac{d^2\widetilde{\psi}}{dz^2}(z_0^*)=\dfrac{2A_2}{(-z_0^*+D_2)^3}
					+\dfrac{2A_3}{(-z_0^*+D_3)^3}-\sum\limits_{i \in S_3}\dfrac{2a_i}{(-z_0^*+d_i)^3}<0.
 \end{align}
 By \eqref{defh2}, \eqref{defh'} and chain rule, we have
    \begin{align}
        &\label{dh=}\dfrac{d\psi}{dz}(z) = \widetilde{\psi}(z)+(z+D_1)\frac{d\widetilde{\psi}}{dz}(z),\;\text{and}\\
					&\label{d2h=}\dfrac{d^2\psi}{dz^2}(z)= \frac{d\widetilde{\psi}}{dz}(z)+ (z+D_1)\frac{d^2\widetilde{\psi}}{dz^2}(z).
    \end{align}
    So, by \eqref{dH=} and \eqref{d2H<}, we have 
				\begin{align}
	\label{dH=*}&\widetilde{\psi}(z_0^*)+(z_0^*+D_1)\frac{d\widetilde{\psi}}{dz}(z_0^*)=0,\;\text{and}\\
 \label{d2H<*}&\frac{d\widetilde{\psi}}{dz}(z_0^*)+ (z_0^*+D_1)\frac{d^2\widetilde{\psi}}{dz^2}(z_0^*)\leq0.
				\end{align}

	Recall that $z_0^*+D_1>0$. By \eqref{final1} and \eqref{dH=*}, we have \eqref{final2}. By \eqref{final2} and \eqref{d2H<*}, we have
 \eqref{final3}.
 
    (Step 4) In this step, we prove that \eqref{final1}, \eqref{final2}, and \eqref{final3} can not hold concurrently.
				Let $\gamma=\dfrac{A_3}{A_2},\; e = \dfrac{-z_0^*
					+D_3}{-z_0^*+D_2},$ and for any $i \in S_3$, let $
				\gamma_{i}=\dfrac{a_i}{A_2}\;\text{and}\; 
				e_{i} = \dfrac{-z_0^*+d_i}{-z_0^*+D_2}$. By \eqref{G1}, we have $\gamma_i>1$ and $\gamma>\sum\limits_{i\in S_3}\gamma_i$.  By \eqref{G2} and \eqref{M=d2a}, we have $e_i>1$ and $e>0$. Then, we divide the both sides of \eqref{final1} by $\frac{A_2}{-z_0^*+D_2}$, and we get \begin{align}\label{final4}
					1+\dfrac{\gamma}{e}
					\leq\sum\limits_{i\in S_3}\dfrac{\gamma_{i}}{e_{i}}.
				\end{align}
				Similarly, by \eqref{final2} and \eqref{final3}, we have 
				\begin{align}
					&\label{final5}1+\dfrac{\gamma}{e^2}
					\geq\sum\limits_{i\in S_3}\dfrac{\gamma_{i}}{e_{i}^2},\;\text{and}\\
					&\label{final5'}1+\dfrac{\gamma}{e^3}
					\leq\sum\limits_{i\in S_3}\dfrac{\gamma_{i}}{e_{i}^3}.
				\end{align}
				By Lemma \ref{finalineq}, \eqref{final4}, \eqref{final5} and \eqref{final5'} lead to a contradiction.
			\end{enumerate}
		\end{enumerate}
\end{enumerate}
	\subsection{Proof of Theorem \ref{thm:main} \ (c)}
	The proof is similar to the proof of (b). So, we put the details in the  supplementary materials \footnote{https://github.com/65536-1024/one-dim}.	
	
	\subsection{Proof of Theorem \ref{thm:main}\;(d)}
	In this case, only one of the four sets $S_1,\; S_2,\; S_3,\; S_4$ is non-empty. So, $g(z)$ is monotone, and hence, $g(z)=0$ has at most one real solution. By Lemma  \ref{diff di}, $G$ admits no multistability.



\end{document}